\documentclass{article}

\usepackage{arxiv}

\usepackage[utf8]{inputenc} 
\usepackage[T1]{fontenc}    
\usepackage{hyperref}       
\usepackage{url}            
\usepackage{booktabs}       
\usepackage{amsfonts}       
\usepackage{nicefrac}       
\usepackage{microtype}      
\usepackage{lipsum}		
\usepackage{natbib}
\usepackage{doi}

\usepackage{amssymb}
\usepackage{gensymb}
\usepackage{multirow}
\usepackage{color,array,dcolumn}
\usepackage{fixltx2e}
\usepackage{stfloats}
\usepackage{url}
\usepackage{amssymb}
\usepackage{gensymb} 
\usepackage{multirow} 

\usepackage{enumitem}
\usepackage{amsmath}
\usepackage{graphicx}
\usepackage{multirow,bigdelim}
\usepackage{booktabs}
\usepackage{textcomp}
\usepackage{adjustbox}
\usepackage{verbatim}
\usepackage{enumitem}
\usepackage{caption}
\usepackage{subcaption}
\usepackage{mathtools}
\usepackage{tikz}
\usetikzlibrary{graphs, graphs.standard, quotes}
\usepackage{todonotes}
\usepackage{wrapfig}
\usepackage{color}
\usepackage[normalem]{ulem}

\newcommand{\vect}[1]{\boldsymbol{\mathrm{#1}}}

\title{Learning Runge-Kutta integration schemes for ODE simulation and identification}

\author{%
  Said Ouala$^1$, Laurent Debreu$^2$, Bertrand Chapron$^3$ \\
  \textbf{Ananda Pascual}$^4$, \textbf{Fabrice Collard}$^5$, \textbf{Lucile Gaultier}$^5$, \textbf{Ronan Fablet}$^1$\\
   $(1)$ IMT Atlantique; Lab-STICC, 29200 Brest, France\\
   $(2)$ AIRSEA INRIA GRA, 38000 Grenoble, France\\ 
   $(3)$ Ifremer, LOPS, 29200 Brest, France\\
   $(4)$ IMEDEA, UIB-CSIC, 07190 Esporles, Spain\\
   $(5)$ ODL, 29200 Brest, France\\
  \texttt{\{said.ouala, ronan.fablet\}@imt-atlantique.fr}\\
  \texttt{Laurent.Debreu@inria.fr},
  \texttt{Bertrand.Chapron@ifremer.fr},\\ \texttt{ananda.pascual@imedea.uib-csic.es} \\
  \texttt{\{dr.fab, lucile.gaultier\}@oceandatalab.com}
}

\renewcommand{\shorttitle}{Learning Runge-Kutta integration schemes for ODE simulation and identification}

\hypersetup{
pdftitle={Learning Runge-Kutta integration schemes for ODE simulation and identification},
pdfsubject={math.NA, stat.ML},
pdfauthor={Said Ouala, Laurent Debreu, Ananda Pascual, Bertrand Chapron,Fabrice Collard, Lucile Gaultier, Ronan Fablet},
pdfkeywords={Dynamical systems, Numerical integration, Automatic differentiation, Runge-Kutta, Data-driven models, Neural ODE},
}
\begin{document}
\maketitle

\begin{abstract}
Deriving analytical solutions of ordinary differential equations is usually restricted to a small subset of problems and numerical techniques are considered. Inevitably, a numerical simulation of a differential equation will then always be distinct from a true analytical solution. An efficient integration scheme shall further not only provide a trajectory throughout a given state, but also be derived to ensure the generated simulation to be close to the analytical one. Consequently, several integration schemes were developed for different classes of differential equations. Unfortunately, when considering the integration of complex non-linear systems, as well as the identification of non-linear equations from data, this choice of the integration scheme is often far from being trivial. In this paper, we propose a novel framework to learn integration schemes that minimize an integration-related cost function. We demonstrate the relevance of the proposed learning-based approach for non-linear equations and include a quantitative analysis w.r.t. classical state-of-the-art integration techniques, especially where the latter may not apply.
\end{abstract}


\section{Introduction}
\label{sec:introduction}

Only linear differential equations and a restricted class of non-linear ones have closed-form analytical solutions. Therefore, for most real-world applications, the resolution of differential equations relies on  
numerical integration algorithms to derive an approximation of the true solution of a given equation.

From a mathematical point of view, applying a numerical integration scheme to a differential equation corresponds to map a continuous time-varying equation to a discrete one \cite{Integrating_traj_book_chapter}. Concequently,  the properties of the simulated solution will depend on the properties of both the continuous equation and the  integration scheme. It is however exteremely important for the properties of the numerical approximation (typically its asymptotic behaviour) to match the true solution according to some  given error and stability criteria.

A large literature discusses such error and stability criteria \cite{Thomas1995} to assess the applicability of an integration scheme for the numerical resolution of a differential equation \cite{friedmann1977efficient,BERLAND20061459,ruuth2006global} with respect to a given discretization and an error tolerance. However, these criteria are most of the time built on linear equations leading to little to no guarantees when considering their behaviour on non-linear systems \cite{Integrating_traj_book_chapter}. In this situation, the integration error depends non-linearly on the state of the system in the phase space which makes the prior characterization of the dynamics not as straightforward as in the linear case (since it becomes equation dependent). Some works developed stability criteria for specific classes of non-linear equations (e.g., B-stability for instance in the case of monotonic equations \cite{B_stability_Book_chap}), however, these criteria are usually restricted to a limited class of equations. 
As a consequence, integration methods are, in practice, often used as black-box routines that simulate a given equation and achieve a given performance criterion. 

From this point of view, new black-boxes, written in the language of data-driven models have quickly joined the arsenal of numerical integration techniques both for the resolution of Ordinary Differential Equations (ODEs) and Partial Differential Equations (PDEs). The success of such techniques is mainly due to their ability to match sufficiently well a performance criterion given a sufficient amount of training data. The main investigations of such approaches relate to the exploration of neural solvers for the resolution of ODEs and PDEs \cite{DBLP:journals/corr/abs-2009-03730,khoo_lu_ying_2020,li2020fourier} and the derivation of some kind of an optimal spatial discretization for the resolution of PDEs \cite{suzuki2019neural,bar2019learning}. Unfortunately, little to no work have been investigated when considering time integration techniques and their role and potential improvements when considering the implementation, resolution and identification of differential equations in a data-driven fashion.

Here, we address the design of time integration schemes as a machine learning problem. We focus on Ordinary Differential Equations (ODEs) and show that one can train numerically data-driven or application-specific integration schemes, possibly jointly to data-driven representations in identification scenarios \cite{wang1998runge,brunton_discovering_2016,fablet_bilinear_2017}. Focusing on an explicit Runge-Kutta formulation, the proposed approach relies on the exploitation of the automatic differentiation \cite{bucker2006automatic} of the computational graphs associated with integration schemes. We demonstrate its relevance for the integration and data-driven identification of non-linear dynamics and assess the properties of the trained integration schemes. One of the main findings of this work is the ability of such methodology to output application-adapted integration schemes without resorting to a complex analysis of the underlying equations. Overall, our key contributions comprise: 
\begin{itemize}
    \item the definition of a generic workflow for integration-related problems where the integration scheme is defined as a solution of an optimization problem;
    \item the data-driven identification of computationally-efficient integration schemes for the time integration of non-linear differential equations;
    \item the data-driven identification of governing laws stated as the joint learning of an approximate ODE and of the associated integration scheme; 
    \item the evaluation of the key features of the trained integration schemes w.r.t. classical state-of-the-art integration techniques, especially through a study of precision, stability and generalization properties for configurations that were never seen during the training phase.
\end{itemize}

The paper is organized as follows. In Section \ref{sec:3SOTA}, we briefly review state-of-the-art integration techniques and the associated performance criteria. Section \ref{sec:ADRK} presents the proposed framework, followed by the experiments and results in Section \ref{sec:Experiments_Linear} and \ref{sec:Exp_L63}. We close this paper with conclusions and perspectives for future works in Section \ref{sec:3Conclusion}.

\section{Integration methods and performance criteria}
\label{sec:3SOTA}
Let us assume a continuous $s$-dimensional dynamical system $\vect{z}_t$ governed by the following non-autonomous time varying ODE
\begin{equation}
     \Dot{\vect{z}}_{t} = f(t,\vect{z}_{t})
     \label{EQ:sec3_ODE_eq}
\end{equation}
Assuming that, given an initial condition $\vect{z}_{t_0}$, we aim to solve this equation for an interval $t\in [t_0,t_f]$, the corresponding solution (or flow) can be written as
\begin{equation}
     \Phi_t(\vect{z}_{t_0}) = \vect{z}_{t_0} + \int_{t_0}^{t}f(w,\vect{z}_{w})dw
     \label{EQ:sec3_Flow_ODE_eq}
\end{equation}
As stated in the introduction, solving the flow integral $\int_{t_0}^{t}f(w,\vect{z}_{w})dw$ is only possible for a small subset of non-linear ODEs. Simulating the differential equation is then done using numerical integration methods. Formally, the interval $t\in [t_0,t_f]$ is discretized using a time-step $h>0$ as $h = \frac{t_f-t_0}{N}$ and $t_n = t_0+n h$, where $0\le n \le N$ an integer and $N$ is the number of grid points, the problem is then formulated as the approximation of the values of variable $\vect{z}_{t}$ at each grid point given the initial condition $\vect{z}_0$ as follows:
\begin{equation}
\left\{
\begin{aligned}
&{\hat{\vect{z}}}_{t_0} = {{\vect{z}}}_{t_0} = {\vect{z}}_{0}\\
&\hat{\vect{z}}_{t_{n+1}} = \Phi_{\mathcal{D},t_{n+1}}(\hat{\vect{z}}_{t_n}) \approx \vect{z}_{t_{n+1}} = \Phi_{t_{n+1}}(\vect{z}_{t_n})
\label{eq:Approx_problem}
\end{aligned}
\right.
\end{equation}
with $\hat{\vect{z}}$ a numerical solution computed using the approximation $\Phi_{\mathcal{D},t_n}$ of the analytical solution (\ref{EQ:sec3_Flow_ODE_eq}) and $\mathcal{D}$ is a given integration scheme.


In this section we will briefly introduce state-of-the-art integration techniques and the associated performance and stability criteria. We focus on single-step explicit schemes as they are one of the most used techniques when considering non-stiff problems.

\subsection{Single-step numerical integration}

Single-step integration schemes can be roughly divided into two main categories: explicit and implicit integration schemes \cite{Integrating_traj_book_chapter}.  Implicit schemes allow larger stability intervals than their explicit counterparts \cite{zhang2016high} which make them extremely valuable when considering stiff problems \cite{Wanner1977,pareschi2000implicit}. By contrast, explicit integration schemes are straightforward to implement and are generally applied  whenever stiffness is not detected within a system.

Single-step explicit integration schemes use a single evaluation of the state of the system at a given grid point to compute the approximation of the solution at the next grid point. Their general form  can be written as follows:
\begin{equation}
\label{EQ:C3_Sing_Step_Exp_Algo}
\hat{\vect{z}}_{t_{n+1}} = \Phi_{\mathcal{D},t_{n+1}}(\hat{\vect{z}}_{t_n}) = \hat{\vect{z}}_{t_{n}} + h\Psi_{\mathcal{D}}(t_{n},\hat{\vect{z}}_{t_{n}},h)
\end{equation}
with $\Psi_{\mathcal{D}}$ an approximation operator related to a truncation of the Taylor expansion of the true solution. We may remind the exact Taylor expansion of the solution of (\ref{EQ:sec3_Flow_ODE_eq}) given the true state $\vect{z}_{t_{n}}$ as:
\begin{equation}
\label{EQ:sec3_Taylor}
\vect{z}_{t_{n+1}} = \vect{z}_{t_{n}} +  \displaystyle \sum_{k=1}^{p = +\infty} h^k\frac{1}{k!}f^{k-1}(t_{n},\vect{z}_{t_{n}})
\end{equation}
with $f^k$ the $k^{th}$-order derivative of ODE operator $f$.

The Forwrad Euler scheme and Runge-Kutta methods are among the most widely used single-step explicit schemes:
\begin{itemize}
    \item \textbf{Forward Euler scheme}: The Forward Euler algorithm may be viewed as a truncation of the Taylor series up to $p = 1$ 
    \begin{equation}
    \label{EQ:sec3_For_Euler}
    \hat{\vect{z}}_{t_{n+1}} = \Phi_{\mathcal{E},t_n}(\hat{\vect{z}}_{t_0}) = \hat{\vect{z}}_{t_{n}} + hf(t_{n},\hat{\vect{z}}^{T}_{t_{n}})
    \end{equation}
    \item \textbf{Explicit Runge-Kutta schemes}: Approximating the integral of (\ref{EQ:sec3_Flow_ODE_eq}) using a low-order truncation of its Taylor series as for instance shown in the Euler algorithm is a quite simple and intuitive approach for approximate solutions of ODEs. However, a simple forward Euler integration can lead to large integration errors and the acceptable integration time-step of such technique is usually small. Considering a high-order truncation usually leads to more accurate and/or stable results. Unfortunately, evaluating high order derivatives of (\ref{EQ:sec3_ODE_eq}) is computationally intensive.
    
    Runge-Kutta methods \cite{RK_solvers}  rely on an iterative evaluation of the ODE function $f$ in order to mimic high order error truncations of the Talyor expansion of a linear equation. A $q$-stage Runge-Kutta ($\mathcal{RK}_q$) technique can be defined by its set of parameters as follows:
\begin{equation}
\label{EQ:RK-q}
\begin{split}
\centering
\mathcal{RK}_q &= \{A = [a_{i,j}]\in \mathbb{R}^{q\times q}, b = [b_i]\in \mathbb{R}^{q}, c = [c_i] \in \mathbb{R}^{q}\\ 
     &\mbox{with }\sum_{j=1}^q a_{i,j} = c_i, 0<c_i<1\mbox{ and  }\sum_{i=1}^q b_{i} = 1 \}
\end{split}
\end{equation}

Given a Runge-Kutta scheme $\mathcal{RK}_q$, the solution of problem (\ref{EQ:sec3_ODE_eq}) can be written as follows:
\begin{equation}
\label{EQ:Sing_Step_Exp_RK}
\hat{\vect{z}}_{t_{n+1}} = \Phi_{\mathcal{RK}_q,t_{n+1}}(\hat{\vect{z}}_{t_{n}})= \hat{\vect{z}}_{t_{n}} + \displaystyle \sum_{i=1}^q b_i k_i
\end{equation}
where $\Phi_{\mathcal{RK}_q,t_{n+1}}(\vect{z}_{t_{n}})$ is the approximate version of (\ref{EQ:sec3_Flow_ODE_eq}) written through the Runge-kutta scheme $\mathcal{RK}_q$. The coefficients $k_i$  are given by:

    \begin{equation}
    \label{EQ:ki_computation}
    k_i = f(t_{n}+c_i h,\hat{\vect{z}}_{t_{n}} + h(\sum_{j=1}^{i-1} a_{i,j} k_j ))
    \end{equation}
    When $q=1$, 
    it simply corresponds to the explicit Euler method.
    For a given number of stages $q$, the coefficients of the Runge-Kutta method  shall satisfy some extra conditions (by matching it to the corresponding Taylor series) to reach a given order $p$ \cite{butcher_1963,wanner1996solving}. Formally, the order $p$ of the Runge-Kutta method  is always inferior or equal to the number of stages $q$. For $q=4$, we can retrieve the well-known Runge-Kutta-4 method, when $p>4$, we need more integration stages $q$ to truly reach a given error order $p$ \cite{wanner1996solving}.
\end{itemize}


Explicit integration schemes approximate the solution at ${t_{n+1}}$ using a forward Taylor expansion, resulting in an expression of $\hat{\vect{z}}_{t_{n+1}}$ as a function of the previous steps (as illustrated for instance by (\ref{EQ:C3_Sing_Step_Exp_Algo})). Implicit schemes \cite{Implicit_RK} in the other hand exploit a backward formulation, resulting in the following general formulation
\begin{equation}
\label{EQ:Sing_Step_Imp_Algo}
\hat{\vect{z}}_{t_{n+1}} = \hat{\vect{z}}_{t_{n}} + h\Psi_{\mathcal{D}}(t_{n+1},\hat{\vect{z}}_{t_{n}},\hat{\vect{z}}_{t_{n+1}},h)
\end{equation}
Note that now, the right hand side of (\ref{EQ:Sing_Step_Imp_Algo}) also depends on $\hat{\vect{z}}_{t_{n+1}}$. The solution of this equation is not as straightforward as in the explicit case since, when considering a non-linear ODE, (\ref{EQ:Sing_Step_Imp_Algo}) becomes non-linear and solving for $\hat{\vect{z}}_{t_{n+1}}$ requires using a non-linear (or linearized) solver \cite{cooper1983iteration}. For a Partial Differential Equation (PDE), it also requires the inversion of a potentially large system. Here, we focus on explicit schemes. We further discuss in Section \ref{sec:3Conclusion} possible applications to explicit schemes. 

\subsection{Integration Errors}
\label{sec:integration_errors}
Integration errors can be seen as the difference between the true analytical solution and the output of the discrete equation at the same time-step (or steps, depending on whether we are considering local or global errors). We may distinguish round-off errors and truncation errors.

The Round-off error is simply the error resulting from using finite-precision arithmetics. This error is inversely proportional to the integration step as a larger integration step would mean a smaller number of time-steps to be computed and stored \cite{Integrating_traj_book_chapter}.

The local truncation error is the error associated with the numerical integration scheme on a single integration time-step \footnote{Please note that the truncation error is derived assuming perfect knowledge of the initial conditions along the temporal grid.} \cite{Integrating_traj_book_chapter}. This error streams from the Taylor expansion as follows. Given the true state at time $t_n$, a $p-$order integration algorithm \footnote{The order of an integration algorithm is determined by matching its Taylor expansion to the Taylor expansion of the true solution of a linear equation. The value of $p$ up to which these two series match and they differ afterwards defines the order of the integration scheme. We may note that this definition implicitly involves a matching criterion which might involve some tolerance or decision threshold.} leads, from the Taylor expansion of the solution at $t_{n+1}$ given by (\ref{EQ:sec3_Taylor}), to the following integration error, referred to as the truncation error $\epsilon_n$ \cite{Integrating_traj_book_chapter}
\begin{equation}
\label{EQ:sec3_Truncation_error}
\epsilon_n = E(p,f,\vect{z}_{t_{n}},t_{n})h^{p+1}
\end{equation}
where $E \in \mathbb{R}$ depends on $p,f,\vect{z}$ and $t_{n}$. The global truncation error is the summation of the local truncation errors for the N grid points. Supposing that the initial condition of each integration step is perfect, it can be computed as $\epsilon^g = \sum_{n=0}^N\epsilon_n$.

\begin{figure}[htb!]
\centering
\centering
 \includegraphics[clip,width=0.5\columnwidth,height=9cm]{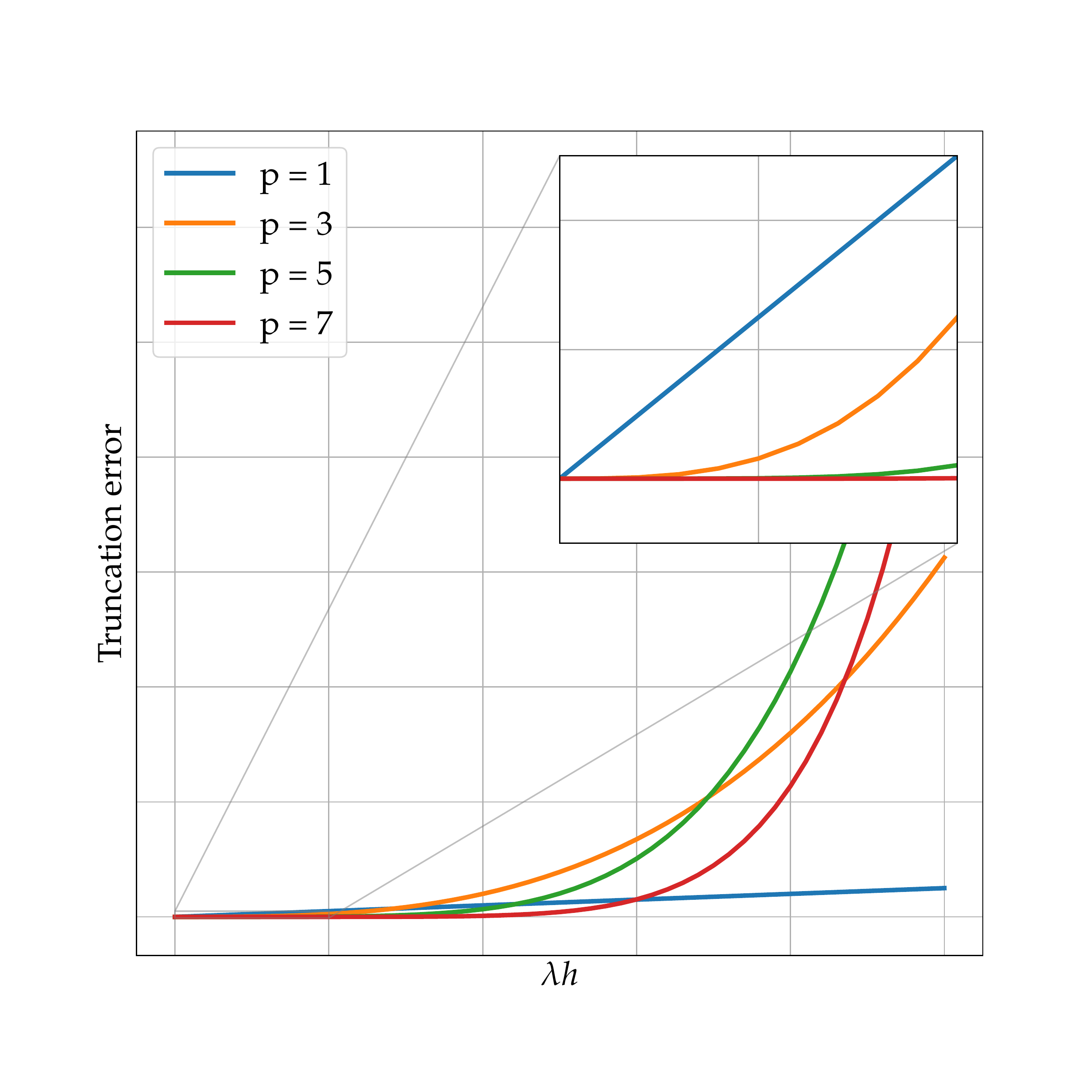}
\caption{{{{ \bf  \em Truncation error of single-step integration schemes.} The truncation error of an integration scheme when integrating the linear equation $\Dot{\vect{z}}_t = \lambda \vect{z}_t$. This error was computed using (\ref{EQ:sec3_Truncation_error_Lin}) for several orders of single-step integration schemes.}
}}
\label{fig:C3_Truncation_Error_Lin}
\end{figure}

We may discuss the truncation error equation. First of all, the truncation error is, for a given numerical integration scheme, proportional to the integration time-step $h$ at a given order. This means that any integration algorithm will produce a larger error as $h$ grows. Furthermore, if the integration time-step is arbitrarily small, high-order algorithms (larger $p$) will achieve a smaller truncation error than low-order algorithms. It is important to note that this statement is only true for small-enough integration time-steps as when $h$ is big enough high-order methods will have higher integration errors due to the term $h^{p+1}$. These considerations can be illustrated on a simple linear ODE. Let us assume that $\vect{z} \in \mathbb{R}$ and $f(t_{n},\vect{z}_{t_{n}}) = \lambda\vect{z}_{t_{n}}$ with $\lambda < 0$. The solution of this dissipative equation is a stable equilibrium point at zero and the corresponding truncation error of a given $p$-order numerical integration scheme becomes
\begin{equation}
\label{EQ:sec3_Truncation_error_Lin}
\epsilon_t = \frac{(\lambda h)^{p+1}}{(p+1)!}\vect{z}_{t_{n}} + O(h^{p+2})
\end{equation}
Figure \ref{fig:C3_Truncation_Error_Lin} depicts a graphical representation of (\ref{EQ:sec3_Truncation_error_Lin}). As stated above, for an arbitrarily small integration time-step (and similarly for an arbitrarily slow time-scale $\lambda$), the higher the order of the integration scheme, the smaller its truncation error. However, and as counter-intuitive as it may sound, above a certain threshold for the values of $(\lambda h)$, lower-order integration schemes may perform better than higher-order ones. This is typically due to the term $(\lambda h)^{p+1}$. These  considerations suggest that choosing an adapted integration scheme for a given application is not an easy task especially when constrained by the integration time-step. In practice,
complex systems involve different values of $\lambda$, which relate to different characteristic time-scales. 
The choice of $h$ and of the order of the time integration scheme may then be difficult to optimize.
One may also notice that the truncation error equation, in the linear case, depends on the state $\vect{z}_{t_{n}}$. In the non-linear case, the error rather depends on the $p^{th}$ derivative of the vector field $f(t_{n},\vect{z}_{t_{n}})$, which generally requires further equation-dependent analyses.

\subsection{Integration Stability}

The characterization of local and global errors provides a relevant framework to assess about the performance of a numerical integration scheme. However, and since they assume perfect knowledge of the initial condition at each time-step, they cannot be used to assess about the convergence properties of an integration scheme when used to simulate a trajectory of an ODE. In practice, local errors accumulate from a time-step to an other and may become unbounded making the integration scheme unstable.

The analysis of the stability of an integration technique usually relies on the study of the numerical solution issued from its application to a first-order linear equation \begin{equation}
     \Dot{\vect{z}}_{t} = \lambda\vect{z}_{t}, \qquad  \vect{z}_{t_0} = \vect{z}_0
     \label{EQ:C3_ODE_eq_lin_exp}
\end{equation}
with $\vect{z}_{t}, \lambda \in \mathbb{C}$ and Real($\lambda$)$ \le 0$.
Since the exact solution of this linear ODE is not growing in time, the corresponding numerical solution must have the same property. This requires stating the discrete system in the form
\begin{equation}
\left\{
\begin{aligned}
&\hat{\vect{z}}_{t_0} = {\vect{z}}_{t_0} = {\vect{z}}_{0}\\
&\hat{\vect{z}}_{t_{n+1}} = \mathcal{R}_{\mathcal{D}}(\lambda h)\hat{\vect{z}}_{t_{n}}\label{EQ:C3_Gain_eq}
\end{aligned}
\right.
\end{equation}
where $\mathcal{R}_{D}$ represents the stability function of the integration scheme. For a $q$-stage explicit Runge-Kutta scheme, the stability function is a polynomial of order less or equal to the number of stages $q$. It can be written as a function of the Runge-Kutta coefficients $\mathcal{RK}_q$ given in (\ref{EQ:RK-q}) as follows \cite{Hairer1996}:
\begin{equation}
     \mathcal{R}_{\mathcal{RK}}(\lambda h) = \sum_{i=1}^{q} \alpha_i (\lambda h)^i = \mathbf{b}^{T}(I-(\lambda h) A)^{-1} \mathbf{1}
     \label{EQ:C3_Gain_RK_analytical}
\end{equation}
with $\alpha_i$ the coefficients of the polynomial $\mathcal{R}_{\mathcal{RK}}(\lambda h)$ and $\mathbf{1}$ a vector of ones with an appropriate dimension.

The region of stability of $\Phi_{\mathcal{D}}$ is given by the values $\lambda h$ for which $|\mathcal{R}_{\mathcal{D}}|\le 1$. We may point out that this widely used stability analysis is only valid for linear equations and does not guarantee the applicability of a given integration scheme for a non-linear ODE.

\section{Automatic Differentiation based Runge-Kutta (ADRK)}
\label{sec:ADRK}

\begin{figure*}
\centering
\includegraphics[width=1.0\textwidth,height=12.0cm]{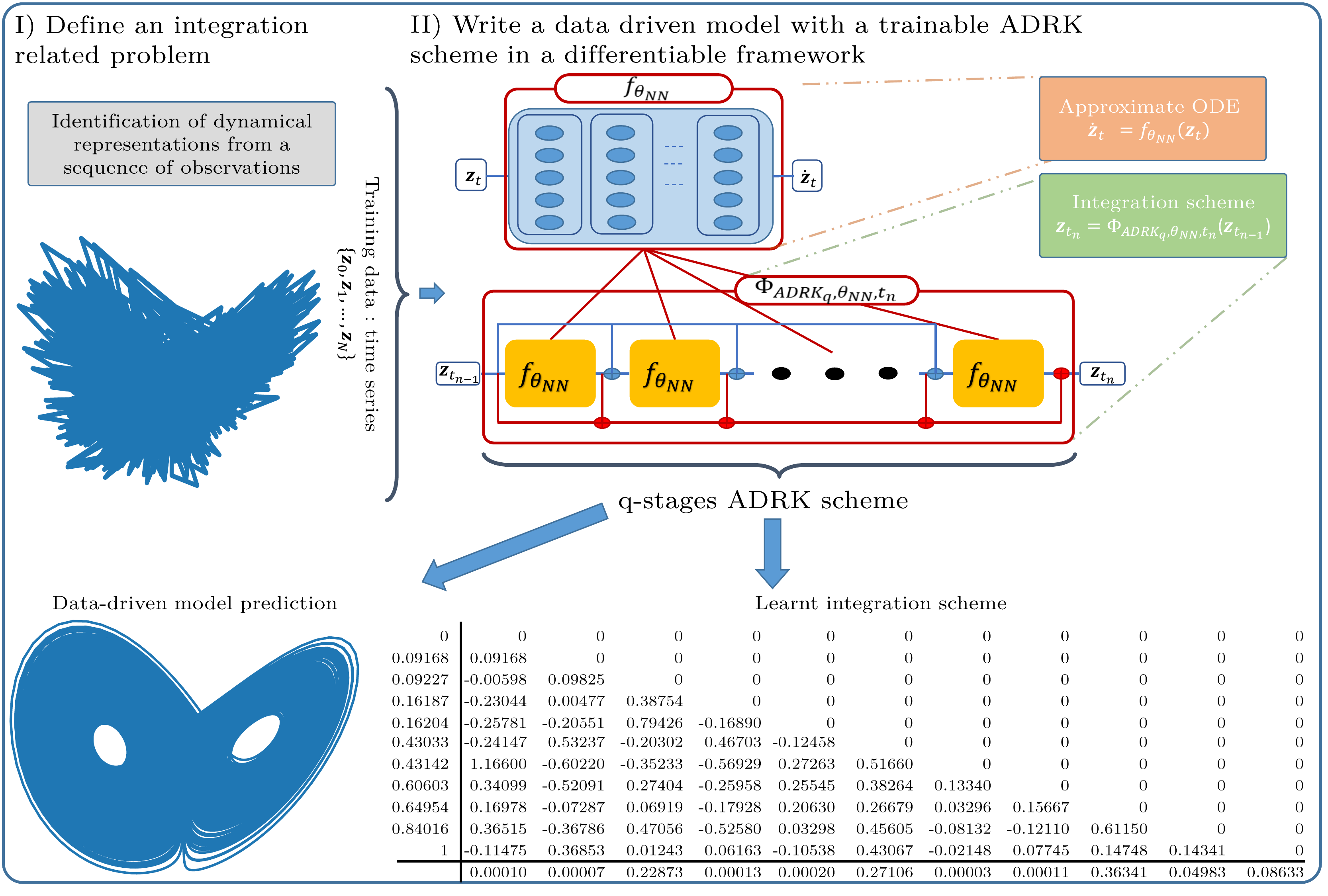}
\caption{{{ \bf   \em Sketch of the proposed ADRK framework.} We highlight the main steps of the proposed framework. We first start by defining an integration related problem, the solution of this problem is written as a Runge-Kutta scheme in a differentiable framework. The parameters of the optimization problem are optimized with respect to a predefined cost. After optimization, the learnt integration scheme can be evaluated using classical stability/precision criterions.}}
\label{fig:Sketch_ADRK}
\end{figure*}

In this section, we introduce a trainable Runge-Kutta integration framework, referred to as Automatic-Differentiation-based Runge-Kutta (ADRK) scheme, to solve integration-related problems
using learning strategies. We first introduce the proposed learning-based framework (Section \ref{sec:learningADRK}). We then address three different applications: the identification of stability-constrained integration schemes (Section \ref{sec:ADRK_optimal_stability}), the integration of ODEs (Section \ref{sec:ADRK_Non-linear_ODE_integration}) and the data-driven identification of ODEs (Section \ref{sec:ADRK_Identification}). 

\subsection{Problem statement}
\label{sec:learningADRK}
In this section, we introduce trainable Runge-Kutta integration schemes framework to solve initial value problems (IVPs). 

Let us recall the general form of an IVP governed by a non-autonomous ODE:
\begin{equation}
\left\{
\begin{aligned}
&\Dot{{\vect{z}}}_{t} = f({t,\vect{z}}_{t})\\
&{\vect{z}}_{t_0} = {\vect{z}}_{0}
\label{eq:C3_ODE_for_RK}
\end{aligned}
\right.
\end{equation}
with ${\vect{z}}_{0}$ is a given initial condition. We may also remind the analytical form of the solution, for an interval $t\in [t_0,t_f]$, as:
\begin{equation}
     \Phi_t(\vect{z}_{t_0}) = \vect{z}_{t_0} + \int_{t_0}^{t}f(w,\vect{z}_{w})dw
     \label{EQ:sec3_Flow_ODE_eq_for_RK}
\end{equation}
As stated in Section \ref{sec:3SOTA}, Runge-Kutta schemes are among the most popular techniques for the numerical resolution of (\ref{EQ:sec3_Flow_ODE_eq_for_RK}).

Given a number of stages $q$, the identification of the parameters of a Runge-Kutta scheme $\mathcal{RK}_q$ given in (\ref{EQ:RK-q}) may derive from different approaches \cite{butcher1996runge,owren1992derivation}. Methods that explicitly enforce stability and/or precision ({\em i.e.} order) constraints on the integration of a linear equation are the most popular ones \cite{ketcheson2013optimal}. They rely on the polynomial form of the integration scheme (given for instance by (\ref{EQ:C3_Gain_RK_analytical})). Such a formulation provides little to no knowledge regarding the stability and/or the performance of the integration scheme on non-linear ODEs. Furthermore, it makes non trivial the choice of the integration scheme when considering identification applications, for which the ODE is not known. 

Here, we adopt a machine learning perspective so that we identify the parameters of a $q$-stage Runge-Kutta scheme which lead to the minimization of a predefined performance criterion. We depict in Figure 
\ref{fig:Sketch_ADRK} a general sketch of the proposed framework, which involves the following key steps:
\begin{itemize}
\item We first state an integration-related problem using a trainable Runge-Kutta scheme with a predefined number of stages. This integration scheme is implemented in a computational framework that supports automatic differentiation \cite{bucker2006automatic} such as pytorch and TensorFlow. Formally, a trainable Runge-Kutta (ADRK) method is defined similarly to (\ref{EQ:RK-q}) by its Runge-Kutta coefficients as follows
\begin{equation}
\label{EQ:ADRK-q}
\begin{split}
\centering
\mathcal{ADRK}_q &= \{A_{\mathcal{ADRK}_q} = [a_{i,j}]\in \mathbb{R}^{q\times q}, b_{\mathcal{ADRK}_q} = [b_i]\in \mathbb{R}^{q}, c_{\mathcal{ADRK}_q} = [c_i] \in \mathbb{R}^{q} \\
&\mbox{ with }\sum_{j=1}^q a_{i,j} = c_i, 0<c_i<1 \mbox{ and  }\sum_{i=1}^q b_{i} = 1 \}
\end{split}
\end{equation}
\item In addition to the parameters of the trainable Runge-Kutta scheme, the considered problem may involve other trainable parameter sets such as unknown parameters of the differential equation. Let us denote by $\nu$ these additional unknown and trainable parameters and $\Theta = \{ \mathcal{ADRK}_q, \nu \}$ the set of all trainable parameters. 
The derivation of optimal parameters  $\widehat{\Theta}$ comes to the following constrained optimization problem with respect to a problem-dependent  criterion $\mathcal{L}$, referred to as the training loss in a deep learning formulation  
\begin{equation}
\begin{array}{c}
    \displaystyle \widehat{\Theta} = \arg \min _\Theta {\mathcal{L}}\left (\Theta \right )\\
    \mbox{subject to} (\ref{EQ:ADRK-q})
    \end{array}
\end{equation}

We will illustrate different problem-dependent definitions for the training loss in the following sections. We may point out that these definitions may involve some priors, such as stability constraints for Runge-Kutta schemes, as well as training data when considering a data-driven identification issues. 

\item Within the considered differentiable framework, we solve this constrained minimization problem using a Stochastic Gradient Descent (SGD) as follows:
\begin{equation}
\label{EQ:Gradient_d_general}
\begin{split}
\centering
\theta^i_{k+1} &= \theta^i_{k} - \mu_i \frac{\partial \mathcal{L}(\Theta) }{\partial \theta^i} 
\end{split}
\end{equation}
with $\mu_i$ a predefined learning rate and $\theta^i \in \Theta$. The constrained optimization for the parameters of the integration scheme is solved in a hard fashion by clipping the Runge-Kutta weights to exactly satisfy the Runge-Kutta constraints (\ref{EQ:ADRK-q}) after each training step $k$. At each iteration $k$, the SGD step typically applies to a subset of the training dataset, referred to as mini-batch in deep learning frameworks \cite{lecun2015deep}. 
\end{itemize}
Deep learning frameworks also 
make highly flexible the definition of the considered training loss, which in turn provide new tools to derive integration schemes according to a variety of criteria as illustrated hereafter. 

\subsection{Stability-constrained ADRK scheme}
\label{sec:ADRK_optimal_stability}

As first case-study of the proposed ADRK framework, we address the derivation of stable Runge-Kutta schemes \cite{ketcheson2013optimal}. Let us consider the following linear ODE:
\begin{equation}
\left\{
\begin{aligned}
&\Dot{{\vect{z}}}_{t} = \lambda{\vect{z}}_{t}\\
&{\vect{z}}_{t_0} = {\vect{z}}_{0}
\label{eq:proposed_mtd_lin_ode}
\end{aligned}
\right.
\end{equation}
with $\lambda \mbox{, } {\vect{z}}_{0} \in \mathbb{C}$ and Real($\lambda$) $\leq 0$. 

The application of a given $q$-stage ADRK scheme (defined by its Runge-Kutta coefficients $\mathcal{ADRK}_q$) in the numerical integration of the linear ODE (\ref{eq:proposed_mtd_lin_ode}) leads to a polynomial form similarly to (\ref{EQ:C3_Gain_eq}) as follow:
\begin{equation}
\left\{
\begin{aligned}
&\hat{\vect{z}}_{t_0} = {\vect{z}}_{t_0} = {\vect{z}}_{0}\\
&\hat{{\vect{z}}}_{t_{n+1}} = \mathcal{R}_{\mathcal{ADRK}_q}(\lambda h) \hat{{\vect{z}}}_{t_{n}}
\label{eq:lin_ode_Gain_ADRK}
\end{aligned}
\right.
\end{equation}
with $\mathcal{R}_{\mathcal{ADRK}_q}$ expressed similarly to (\ref{EQ:C3_Gain_RK_analytical}) as a function of the coefficients $\mathcal{ADRK}_q$ as follows:
\begin{equation}
     \mathcal{R}_{\mathcal{ADRK}}(\lambda h) = \sum_{i=1}^{q} \alpha_i (\lambda h)^i = \mathbf{b}_{\mathcal{ADRK}}^{T}(I-(\lambda h) A_{\mathcal{ADRK}})^{-1} \mathbf{1}
     \label{EQ:C3_Gain_ADRK_analytical}
\end{equation}
The derivation of the integration coefficients is classically done by i) following section \ref{sec:integration_errors}, matching the coefficients $\alpha_i$ to $\frac{1}{i!}$ up to a given order $\hat{p}$ or ii) by optimizing the coefficients $\alpha_i$ to match a stability region $|\mathcal{R}_{\mathcal{ADRK}_q}(\lambda h)| \leq 1$ for a range of values of $\lambda h$. In the latter case, most case studies in previous works focus on specific cases such as dissipative processes {\em i.e.} $\lambda \in \mathbb{R}$ with $\lambda < 0$ (also stated when considering the derivation of stability polynomials as real axis inclusion \cite{abdulle2000roots,abdulle2002fourth,bogatyrev2005effective,ketcheson2013optimal}, for which the maximal stable region of a q-stage Runge-Kutta method is reached, for first order schemes, for $\lambda h = 2q^2$ \cite{ketcheson2013optimal}) and advection processes characterized by $\lambda \in \mathbb{C}$ with Real($\lambda$) $= 0$ (imaginary axis inclusion \cite{KINNMARK198487,KINNMARK1984181,mead1999optimal,ketcheson2013optimal} for which the maximal stable region of a q stages Runge-Kutta method is reached, for first order schemes, for $\lambda h = q-1$ \cite{ketcheson2013optimal}) and other geometries (disk regions, thin rectangles ...etc.). Overall, such a methodology requires the definition of a potentially complex optimization procedure and usually the choice of a polynomial basis for $\mathcal{R}_{\mathcal{ADRK}_q}$.

We apply the  proposed ADRK framework to derive  stability polynomials without relying on a polynomial basis or on a specific optimization technique. Following the general workflow described in the previous section, the proposed framework relies on  the following steps:
\begin{itemize}
\item We first define a linear equation (\ref{eq:proposed_mtd_lin_ode}) for a specific value of $\lambda$ with Real$(\lambda)\leq 0$. We consider $\lambda \in \mathbb{R}$ when deriving optimal stability polynomial in the real axis and similarly $\lambda \in \mathbb{C}$ with Real$(\lambda) =0$ when considering imaginary axis stability inclusion. Other geometries can be considered.

\item the set of trainable parameters comprises the ADRK coefficients \begin{equation}
\label{EQ:trainable_params_polynomes}
\begin{split}
\centering
\Theta &= \{\mathcal{ADRK}_q \}
\end{split}
\end{equation}

\item As training loss $\mathcal{L}(\Theta)$, we consider the following stability-related criterion, such that the integration scheme will be stable or values of $\lambda h$ up to a maximum step-size $h_f$: 
\begin{equation}
\label{EQ:polynom_stable_crit}
\begin{split}
\centering
\mathcal{L}(\Theta) &= \displaystyle  \max\bigg(0,\sum_{h=0}^{h_f} (|\mathcal{R}_{ADRK}(h\lambda)| - 1) \bigg)
\end{split}
\end{equation}
If the minimization of this training loss reaches 0, it guarantees that the trained ADRK scheme is stable over a finite set of time steps from $0$ to $h_f$. To  enforce the stability over a continuous range of values, we noted experimentally that we could gradually decrease the descritezation spacing of $h$, typically up to $h_f/100$, until the training loss no longer changes at zero.


The maximum step size $h_f$ is defined depending on the integration problem. For example, and following \cite{ketcheson2013optimal}, when considering $\lambda \in \mathbb{R}$ with Real$(\lambda) <0$ {\em i.e.} real axis inclusion, a maximum value can be set as $h_f = 2q^2/|\lambda|$. Similarly, if $\lambda$ is purely imaginary, the limit values of stable step-size can be set as $h_f = (q-1)/|\lambda|$. The value of the stability polynomial $\mathcal{R}_{ADRK}$ can be either computed using the analytical formulation (\ref{EQ:C3_Gain_ADRK_analytical}) or based on a ratio between ${{\vect{z}}}_{t_{n+1}}$ and ${{\vect{z}}}_{t_{n}}$ given in (\ref{eq:lin_ode_Gain_ADRK}). Both these two formulations lead to the same value of $\mathcal{R}_{ADRK}$, however, the analytical formulation (\ref{EQ:C3_Gain_ADRK_analytical}), of the stability polynomial, have a faster converges rate when considering large values of $\lambda h$. In this work, real axis inclusion case studies were implemented using the analytical formulation (\ref{EQ:C3_Gain_ADRK_analytical}) and imaginary axis inclusion experiments were implemented based on a ratio between ${{\vect{z}}}_{t_{n+1}}$ and ${{\vect{z}}}_{t_{n}}$.

\end{itemize}

The above framework leads to the derivation of stable Runge-Kutta methods up to a given stability limit. We may apply the same workflow for the derivation of high-order integration schemes. Criterion (\ref{EQ:polynom_stable_crit}) can be replaced by matching the ADRK coefficients $\alpha_i$ to $\frac{1}{i!}$ up to a given order $\hat{p}$. An example of application of this framework is given in section \ref{sec:Experiments_Linear_stability}

\subsection{ODE-adapted ADRK scheme}
\label{sec:ADRK_Non-linear_ODE_integration}
Given a non-linear ODE in (\ref{eq:C3_ODE_for_RK}), the stability and precision properties discussed above provide little to no knowledge regarding the performance of an integration scheme (specifically a Runge-Kutta one) in the time integration of non-linear equations. From this point of view, adaptive step-size solvers are typically used as an efficient framework that deals with non-linear differential equations. This however comes at the expense of an increased computational complexity. 
Here, we apply the proposed ADRK scheme to derive computationally-efficient integration schemes for a given non-linear differential equation.

The proposed ADRK technique can exploit simulations from complex integration tools (adaptive solvers, integration schemes with a high number of stages) to optimize adapted integration routines that run at a smaller cost. Formally, let us assume that we are provided with representative time series of a full state vector $\{\vect{z}\}$ of the ODE (\ref{eq:C3_ODE_for_RK}) with a given time sampling rate $h$. This training dataset may result from the application of a predefined high-complexity solver. For the sake of simplicity, we consider below a single time series of length $N+1$, $\{\vect{z}_{t_0},\vect{z}_{t_1},\ldots,\vect{z}_{t_N}\}$ with $t_n = t_0 + nh$. We could also consider similarly a dataset formed by different time series possibly of varying lengths. 

We apply the proposed ADRK framework to learn a $q$-stage ADRK scheme which shall minimize the forecasting error. This leads to the following workflow:
\begin{itemize}
\item We first implement the non-linear equation (\ref{eq:C3_ODE_for_RK}) in a differentiable framework.  

\item We define the parameters that require optimization as the ADRK coefficients \begin{equation}
\label{EQ:trainable_params_non_lin_integration}
\begin{split}
\centering
\Theta &= \{ \mathcal{ADRK}_q \}
\end{split}
\end{equation}

\item We specify the training loss  $\mathcal{L}(\Theta)$ as a trajectory matching  cost such that the ADRK scheme will minimize the simulation error of a given trajectory $\{\vect{z}_0,\vect{z}_1,\ldots,\vect{z}_N\}$ as follows:
\begin{equation}
\label{EQ:cost_function_integration}
\begin{split}
\centering
\mathcal{L} &= \displaystyle \sum_{n=1}^N \| \vect{z}_{t_{n}} - \Phi_{\mathcal{ADRK}_q,t_n}(\vect{z}_{t_{n-1}})\|  
\end{split}
\end{equation}
with $\Phi_{\mathcal{ADRK}_q,t_n}(\vect{z}_{t_{n-1}})$ is the numerical ADRK integration of (\ref{eq:C3_ODE_for_RK}).
\end{itemize}

Once the training procedure (\ref{EQ:Gradient_d_general}) has converged, we may apply the trained ADRK scheme to compute relevant numerical simulations for the ODE (\ref{eq:C3_ODE_for_RK}). 
The application of this framework is reported in sections \ref{sec:ODE_integration_linear} and \ref{sec:ODE_integration_nlinear} when considering respectively the integration of linear and non-linear ODEs.

\subsection{Joint identification of ADRK schemes and ODE models}
\label{sec:ADRK_Identification}
Deriving an ODE representation that reproduces the variability of a given dataset is a substantially different, and obviously harder, problem than the integration of a known equation. Particularly, choosing an integration algorithm for a known non-linear ODE may be complex 
and selecting an integration technique to learn an ODE appears even more complex. Using adaptive step size solvers is a practical solution of the integration issue, however, adaptive step size solvers as proposed in \cite{chen2018neural} may be subject to memory and stability issues \cite{Acc_rate_grad_NODE1,Acc_rate_grad_NODE2}. 
The proposed framework appears as an interesting alternative to classical integration routines since an ADRK scheme can be learnt jointly with an approximate data-driven model from data.

Similarly to the previous case-study, let us assume that we are provided with a time series of $N+1$ samples of the full state vector $\{\vect{z}_{t_0},\vect{z}_{t_1},\ldots,\vect{z}_{t_N}\}$ with $t_n = t_0 + nh$ and $h$ a given time sampling rate.  Contrary to the previous case, we assume that this state vector is governed by an unknown ODE. We apply the proposed framework to address the joint identification of the governing ODE and of the associated integration scheme. We proceed as follows: 
\begin{itemize}
\item We first define a parameterization of the data-driven ODE \begin{equation}
     \Dot{\vect{z}}_{t} = f_{\theta_{NN}}(t,\vect{z}_{t})
     \label{EQ:sec3_ODE_eq_DD}
\end{equation}
in a differentiable framework, which usually corresponds to the definition of a neural network architecture with parameters $\theta_{NN}$. This archictecture may rely on physics-informed parameterizations \cite{brunton_discovering_2016,raissi_multistep,fablet_bilinear_2017,debezenac_deep_physical} or  on fully data-driven formulations. Regarding the numerical integration of the ODE (\ref{EQ:sec3_ODE_eq_DD}), we aim to apply a trainable  $q$-stage ADRK scheme $\mathcal{ADRK}_q$ given by (\ref{EQ:ADRK-q}). 


\item Here, we aim to jointly learn the parameters of the ADRK integration scheme as well as of the ODE, such that the set of trainable parameters $\Theta$ are given by
\begin{equation}
\label{EQ:trainable_params_DD}
\begin{split}
\centering
\Theta &= \{\mathcal{ADRK}_q, \theta_{NN} \}
\end{split}
\end{equation}

\item Similarly to the previous case-study, the training procedure comes to optimize all trainable parameters with a view to minimizing the forecasting error from a predefined initial condition, which leads to the following training loss:
\begin{equation}
\label{EQ:cost_function_DD}
\begin{split}
\centering
\mathcal{L}(\Theta) &= \displaystyle \sum_{n=1}^N \| \vect{z}_{t_{n}} - \Phi_{\mathcal{ADRK}_q,\theta_{NN},t_n}(\vect{z}_{t_{n-1}})\|  
\end{split}
\end{equation}
with $\Phi_{\mathcal{ADRK}_q,\theta_{NN},t_n}(\vect{z}_{t_{n-1}})$ the numerical ADRK integration of (\ref{EQ:sec3_ODE_eq_DD}).
\end{itemize}
We report an example of application of this framework to the data-driven identification of ODE representations of chaotic dynamics when the available observation data involve a scarce time sampling \ref{sec:identification_L63}

\section{linear ODE case-study}
\label{sec:Experiments_Linear}
In order to illustrate the key principles of the proposed framework, let us consider the following linear ODE:
\begin{equation}
\left\{
\begin{aligned}
&\Dot{{\vect{z}}}_{t} = \lambda{\vect{z}}_{t}\\
&{\vect{z}}_{t_0} = {\vect{z}}_{0}
\label{eq:C3_exp_lin_ode}
\end{aligned}
\right.
\end{equation}
with $\lambda$ either a real or complex constant.

We first perform numerical experiments of the proposed framework to learn stability-constrained  ADRK scheme as introduced in Section \ref{sec:ADRK_optimal_stability}. In a second experiment, we address the optimization of an ADRK scheme for the integration of this linear ODE,  $\lambda$ being known.

\subsection{Learning stability-constrained integration schemes}
\label{sec:Experiments_Linear_stability}

We first report numerical experiments for the learning of stability-constrained ADRK schemes for a given range of integration time-steps. 
Specifically, we investigate in this section the ability of an $\mathcal{ADRK}_q$ scheme in the derivation of a stable polynomial $\mathcal{R}_{\mathcal{ADRK}_q}(\lambda h)$ that can integrate (\ref{eq:C3_exp_lin_ode}) up to the the limit of $\lambda h = 2q^2$ for real axis inclusion ({\em i.e} $\lambda \in \mathbb{R}$) and $\lambda h = q-1$ for imaginary axis inclusion ({\em i.e} $\lambda \in \mathbb{C}$ with Real($\lambda$) = 0). For this purpose, the formulation of section \ref{sec:ADRK_optimal_stability} is studied on a 4-stage ADRK scheme. We include additional experiments, with a higher number of stages, in Supplementaray Material \ref{app:Stability_poly}.

\subsubsection{Real-axis inclusion}
In this experiment, we state a $\mathcal{ADRK}_4$ scheme in the framework described in section \ref{sec:ADRK_optimal_stability} with $h_f =2q^2/|\lambda| = 32/|\lambda|$ and $\lambda = -4$. This case-study relates to the design of stable integration schemes for dissipation processes.
Figure \ref{fig:C3_stability_region_2d} illustrates the stability region of the trained ADRK scheme. A visual comparison of the stability region of this scheme with respect to the classical Euler and Runge Kutta 4 is also given in figure \ref{fig:c3_stabil_region_1D_stab_constrained}. 

The trained $\mathcal{ADRK}_4$ is able to reach the upper bound of $\lambda h = 2q^2$. Furthermore, we can write the polynomial form of the trained ADRK scheme as follows:

\begin{figure}
\centering
\includegraphics[clip,width=0.5\columnwidth,height=9cm]{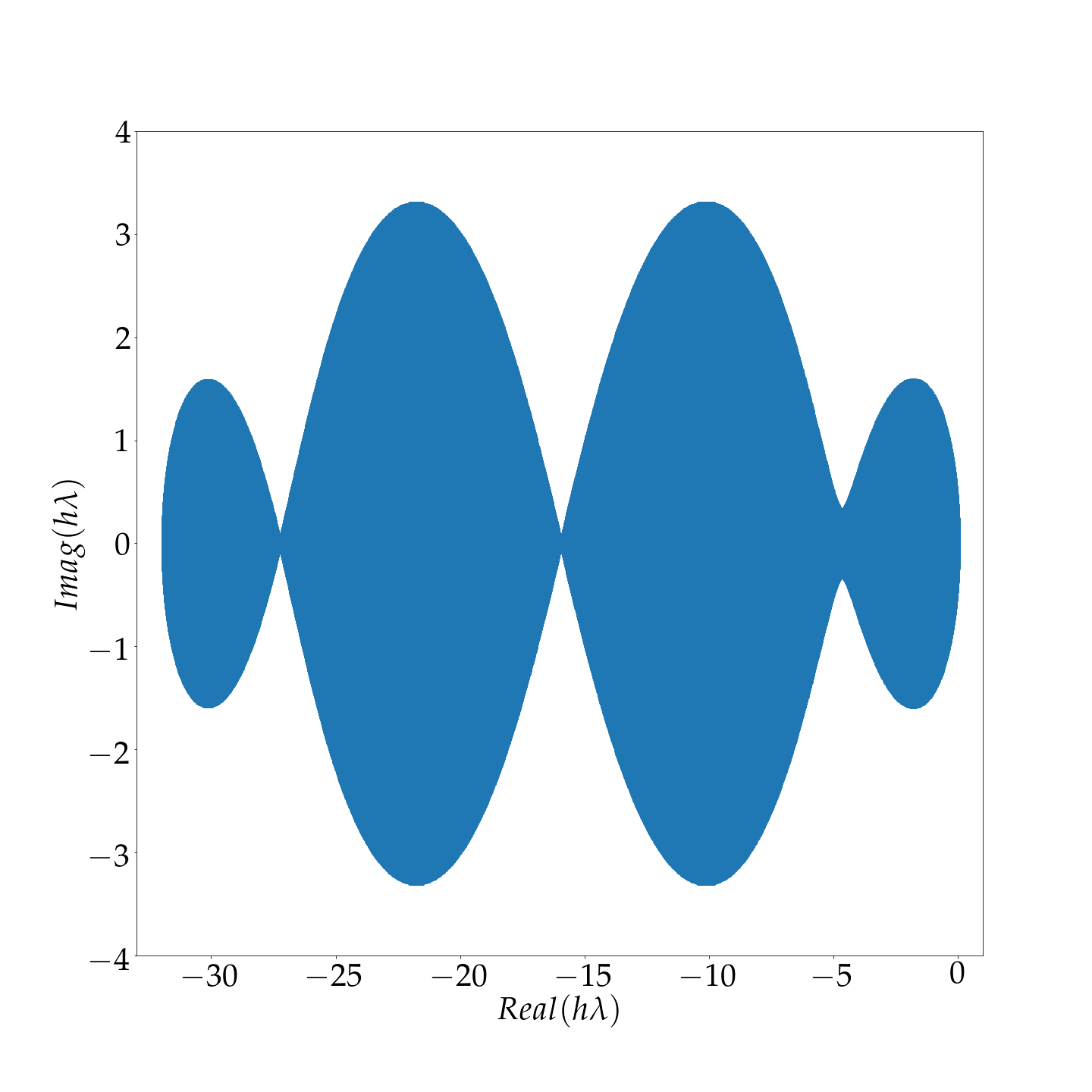}
\caption{{{{ \bf  \em Two dimensional stability region of the proposed $\mathcal{ADRK}_4$ integration scheme for dissipation processes}: The $\mathcal{ADRK}_4$ scheme was trained based on a stability criterion, up to $\lambda h = 2q^2$ for $\lambda \in \mathbb{R}$, $\lambda <0$.}
}}
\label{fig:C3_stability_region_2d}
\end{figure}

\begin{figure}[h]
\centering
  \includegraphics[clip,width=0.7\columnwidth,height=8cm]{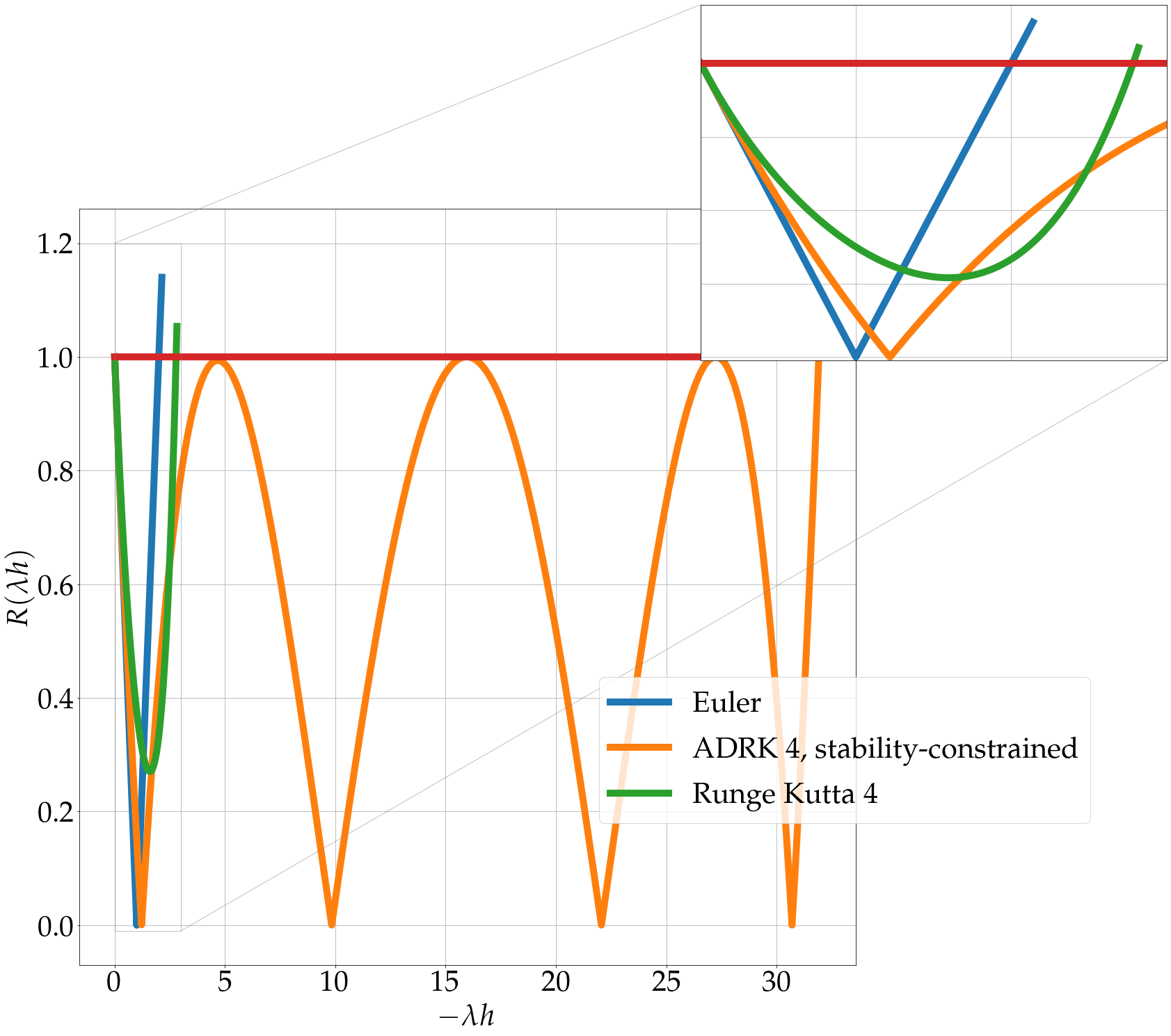}%
\caption{{{{ \bf  \em One-dimensional stability region of the proposed stability constrained ADRK integration scheme with respect to Euler and Runge-Kutta 4 schemes}. We highlight the stability region of the stability constrained ADRK scheme (given in the complex plan in figure \ref{fig:C3_stability_region_2d}) in the real plan with respect to the Euler and Runge-Kutta 4 integration techniques.}
}}
\label{fig:c3_stabil_region_1D_stab_constrained}
\end{figure}

\begin{align}
\begin{split}
\mathcal{R}_{\mathcal{ADRK}_4}(\lambda h)
              =& 1 + \lambda h + 0.1567 (\lambda h)^{2} + 0.0786 (\lambda h)^{3} + 0.0001 (\lambda h)^{4}
\end{split}
\label{EQ:C3_Taylor_gain_RINN_opti_real_axis}
\end{align}
We can remind the Taylor series of the true solution $exp(\lambda h)$ 

\begin{equation}
\label{EQ:C3_Taylor_Exp}
exp(\lambda h) = \displaystyle \sum_{p=0}^{p = \infty} \frac{(\lambda h)^p}{p!}
              = 1 + \lambda h + \frac{1}{2}(\lambda h)^2 + \frac{1}{6}(\lambda h)^3 + ...
\end{equation}
As discussed in Section \ref{sec:integration_errors}, from the comparison of (\ref{EQ:C3_Taylor_gain_RINN_opti_real_axis}) to the expansion of the true solution (\ref{EQ:C3_Taylor_Exp}), this stability-constrained scheme is only first order, which  
can be interpreted as trading the precision performance of the scheme for a larger stability region. 
Related works \cite{ketcheson2013optimal} drew similar observations, as the stability upper bound of $\lambda h = 2q^2$ with $\lambda < 0$ can only be reached using first-order integration schemes.


\subsubsection{Imaginary-axis inclusion}
\begin{figure}[h]
\centering
  \includegraphics[clip,width=0.5\columnwidth,height=9cm]{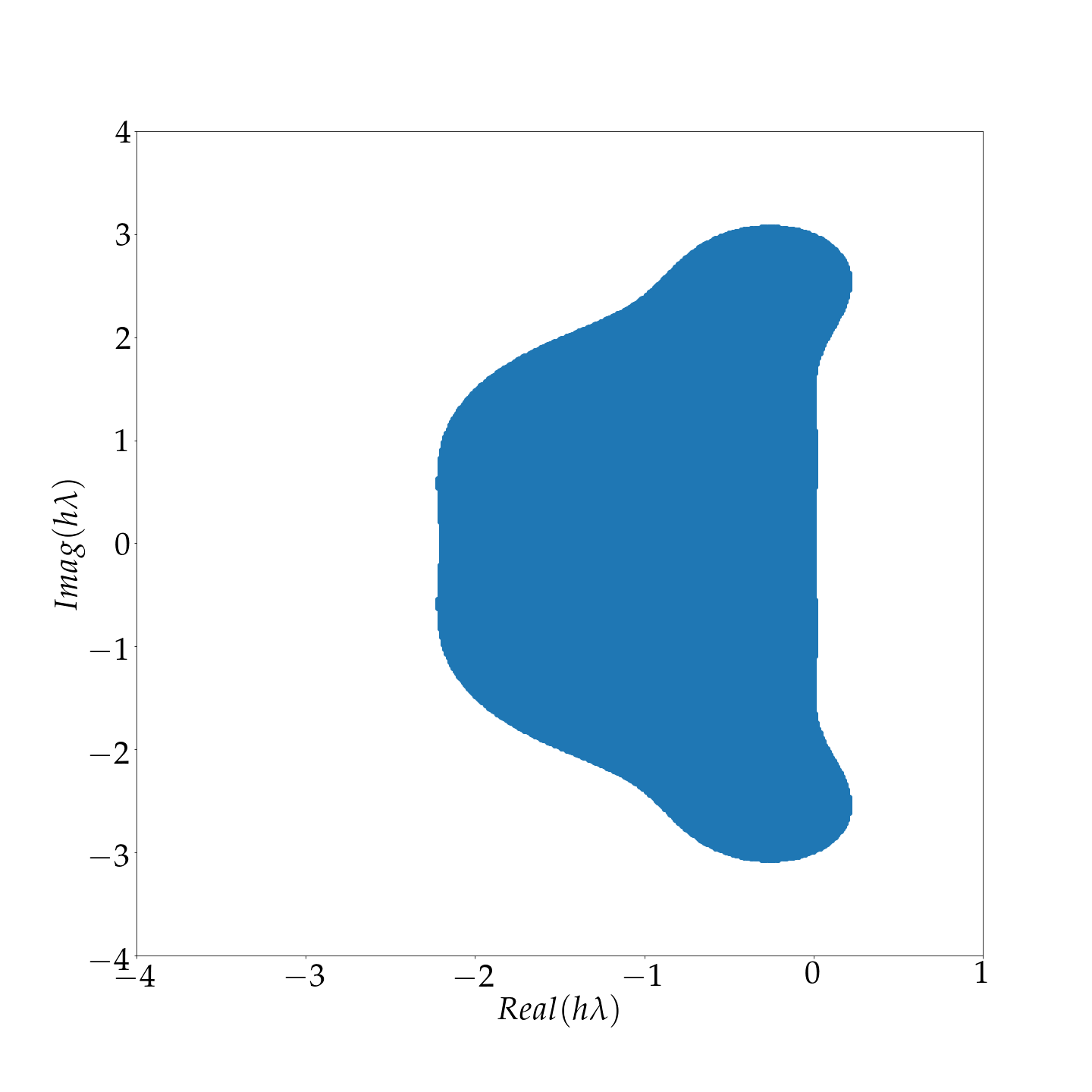}%
\caption{{{{ \bf  \em Two dimensional stability region of the proposed $\mathcal{ADRK}_4$ integration scheme} for advection processes: The $\mathcal{ADRK}_4$ scheme was trained based on a stability criterion, up to $\lambda h = q-1$ for $\lambda \in \mathbb{C}$ with Real($\lambda$)$=0$.}
}}
\label{fig:c3_stability_ADRK4_imaginary_inclusion}
\end{figure}

Here, we apply the stability-constrained $\mathcal{ADRK}_4$ framework described in section \ref{sec:ADRK_optimal_stability} with $h_f = 3/|\lambda|=(q-1)/|\lambda|)$ and $\lambda = -j$ ($j^2 = -1$), which refers to advection process \cite{KINNMARK1984181}. Figure \ref{fig:c3_stability_ADRK4_imaginary_inclusion} illustrates the stability region of the trained ADRK scheme. The trained $\mathcal{ADRK}_4$ is again able to reach the theoretical upper bound of $\lambda h = q-1$ for imaginary-axis inclusion. We can write the polynomial form of the trained ADRK scheme as follows: 
\begin{align}
\begin{split}
\mathcal{R}_{\mathcal{ADRK}_4}(\lambda h)
              =& 1 + \lambda h + 0.5478 (\lambda h)^{2} + 0.1490 (\lambda h)^{3} + 0.0480(\lambda h)^{4}
\end{split}
\label{EQ:C3_Taylor_gain_RINN_opti_imaginary_axis}
\end{align}

\begin{figure*}
\centering
  \subfloat[]{%
\includegraphics[width=0.5\textwidth]{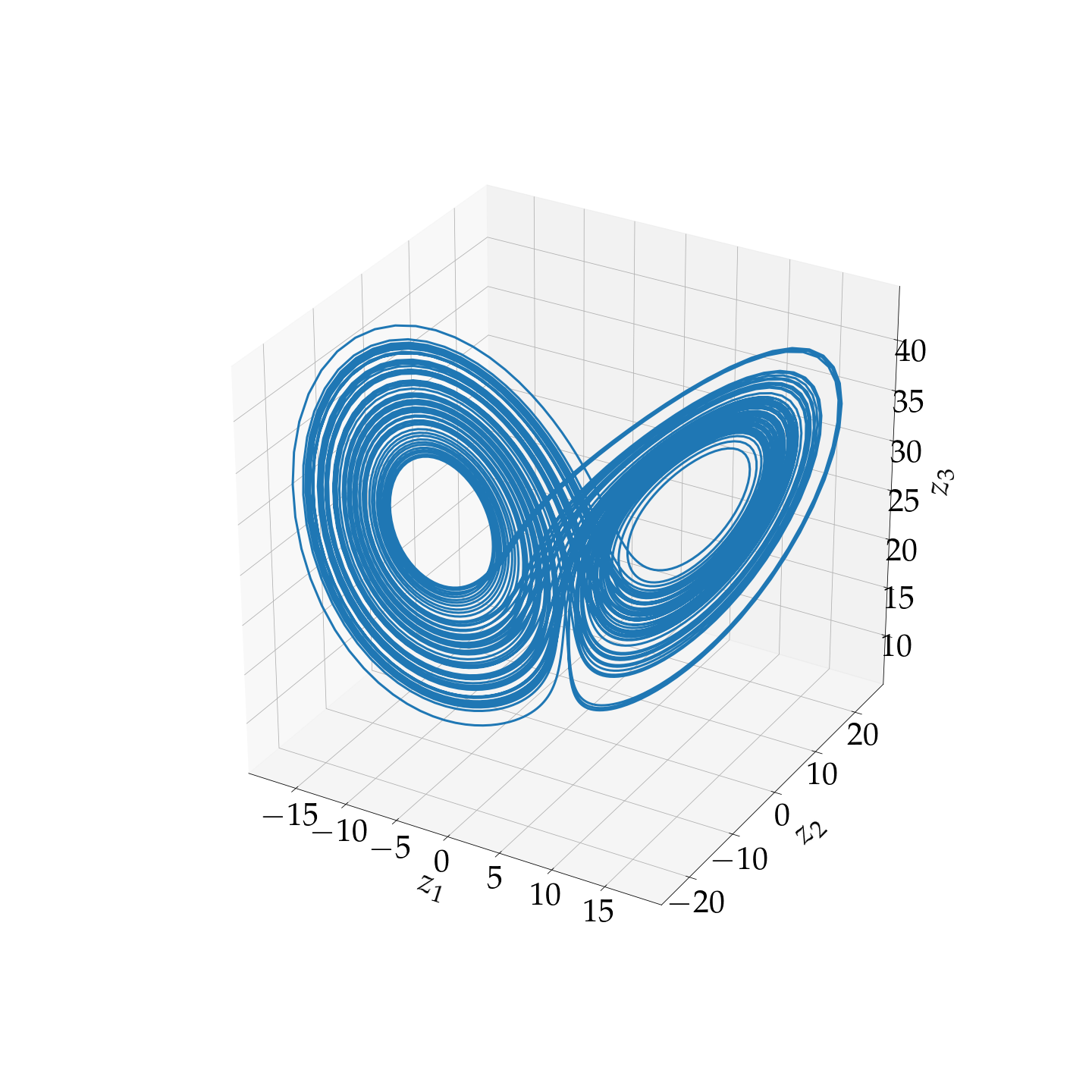}}
  \subfloat[]{%
\includegraphics[width=0.5\textwidth]{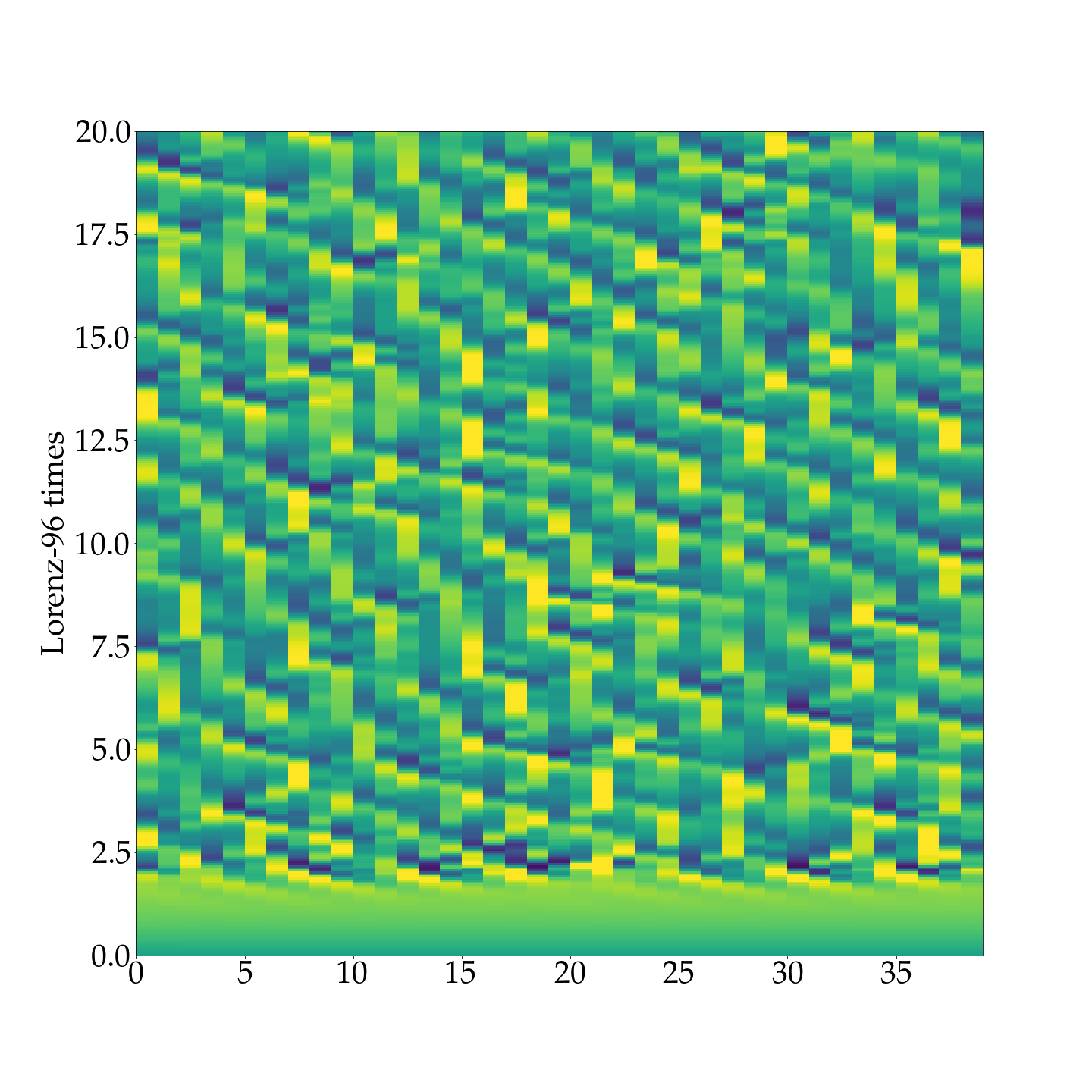}}\\
  \subfloat[]{%
\includegraphics[width=0.5\textwidth]{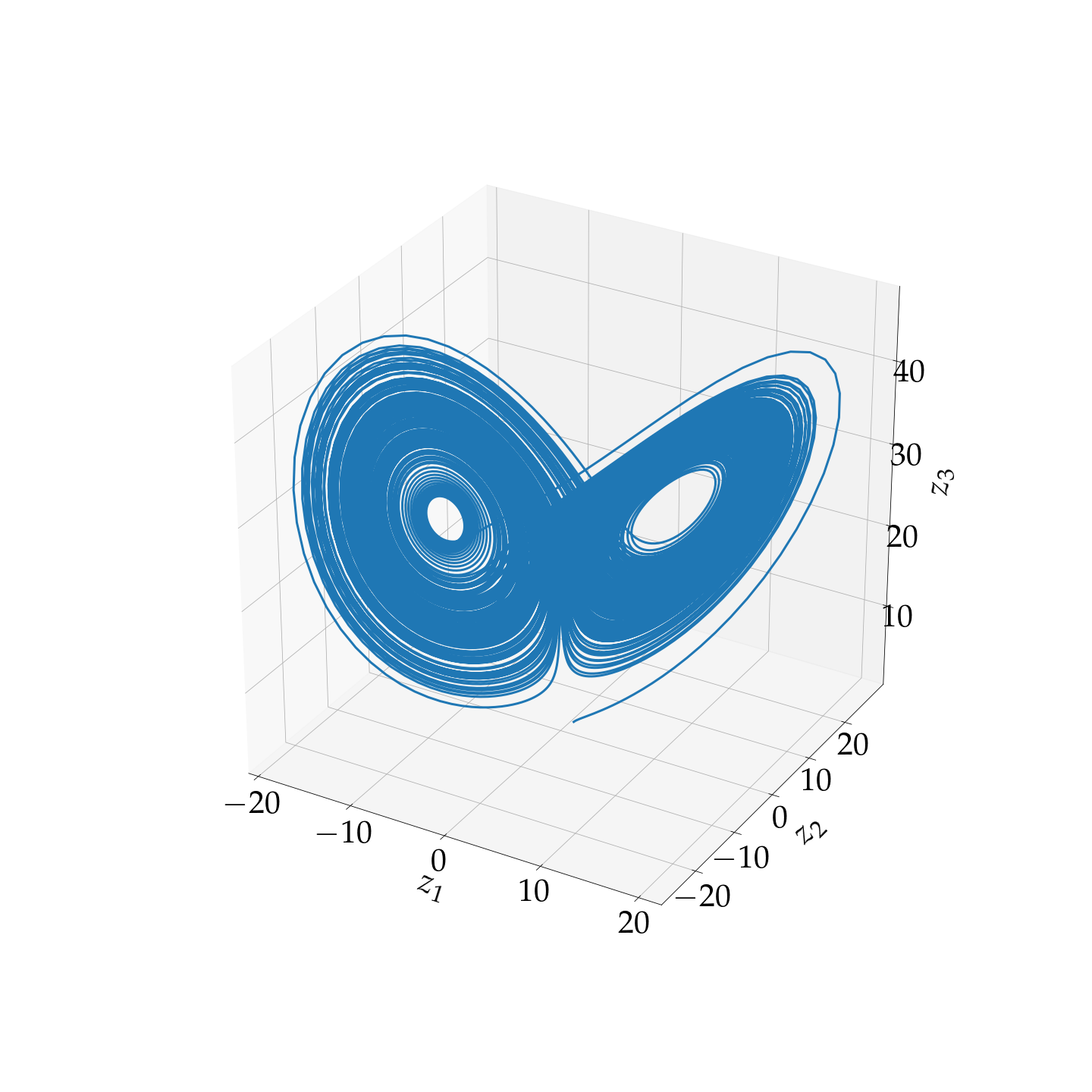}}
  \subfloat[]{%
\includegraphics[width=0.5\textwidth]{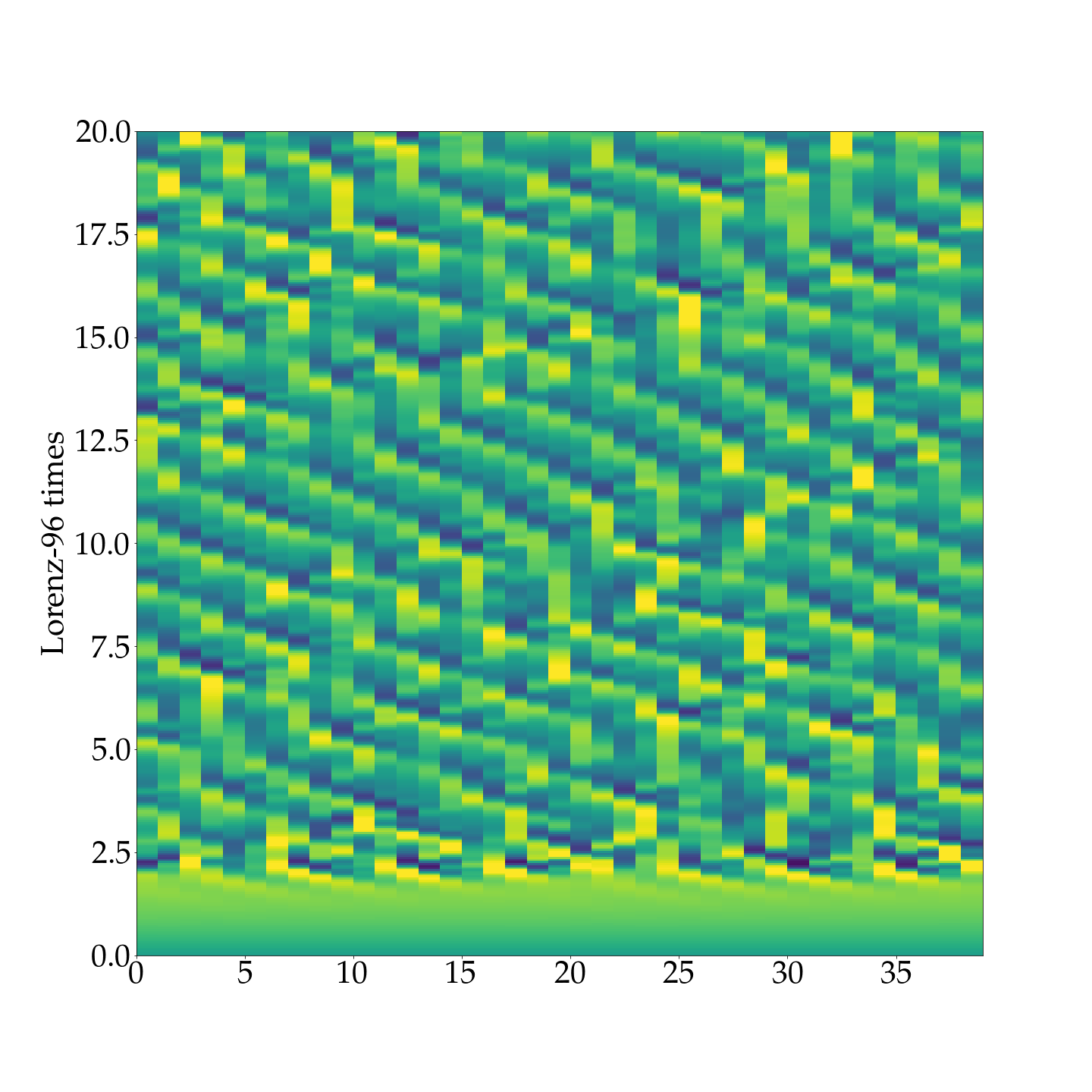}}
\caption{{{{ \bf  \em Application of the stability-constrained $\mathcal{ADRK}_4$ integration scheme on non-linear dynamics}. Integration of Lorenz 63 (a) and 96 (b) equations using our stability-constrained $\mathcal{ADRK}_4$ integration scheme given by the stability polynomial (\ref{EQ:C3_Taylor_gain_RINN_opti_imaginary_axis}). Figures (c) and (d) correspond to the integration of the Lorenz 63 and 96 models respectively using the classical adaptive step-size Adams technique.}
}}
\label{fig:C3_Lin_RINN_On_NL_ODES}
\end{figure*}

The learned integration scheme, represented here for instance by the stability polynomial given in (\ref{EQ:C3_Taylor_gain_RINN_opti_imaginary_axis}), is close to a second order scheme, and shows good stability properties along the real axis (the scheme is still stable slighthly above Real($\lambda h$)$ = -2$) even if it was only trained on purely imaginary values of $\lambda h$. This is principally due to the consistency condition $\sum_{i = 1}^{q} b_{\mathcal{ADRK}_i} = 1$ enforced during the training phase. Finnaly, this stability-constrained integration scheme can also be used to solve other ODEs. Fig. \ref{fig:C3_Lin_RINN_On_NL_ODES} shows the simulated trajectories of Lorenz-63 \cite{lorenz_deterministic_1963} and Lorenz-96 \cite{lorenz1996predictability} systems using the stability-constrained $\mathcal{ADRK}_4$ scheme (given by the stability polynomial (\ref{EQ:C3_Taylor_gain_RINN_opti_imaginary_axis})). Using this integration scheme on such non-linear equations leads to less accurate simulations when compared to the classical fourth order Runge-Kutta since the learnt scheme is only approximately second order and thus less accurate. The learnt scheme is however more stable when considering larger integration time-steps.

\subsection{ADRK scheme optimized for a known value of $\lambda$}
\label{sec:ODE_integration_linear}

Assuming the value of parameter $\lambda$ is known in (\ref{eq:C3_exp_lin_ode}), here $\lambda = -2$, we address the derivation of an optimal 4-stage ADRK to reproduce a simulation of (\ref{eq:C3_exp_lin_ode}). This simulation is computed using the analytical solution of (\ref{eq:C3_exp_lin_ode}) from $t = 0$ up to $t = 10$ with a sampling rate $h = 0.01$. The initial condition $\vect{z}_{0}$ was set to $0.5$.  

\begin{figure}[h]
\centering
  \includegraphics[clip,width=0.5\columnwidth,height=9cm]{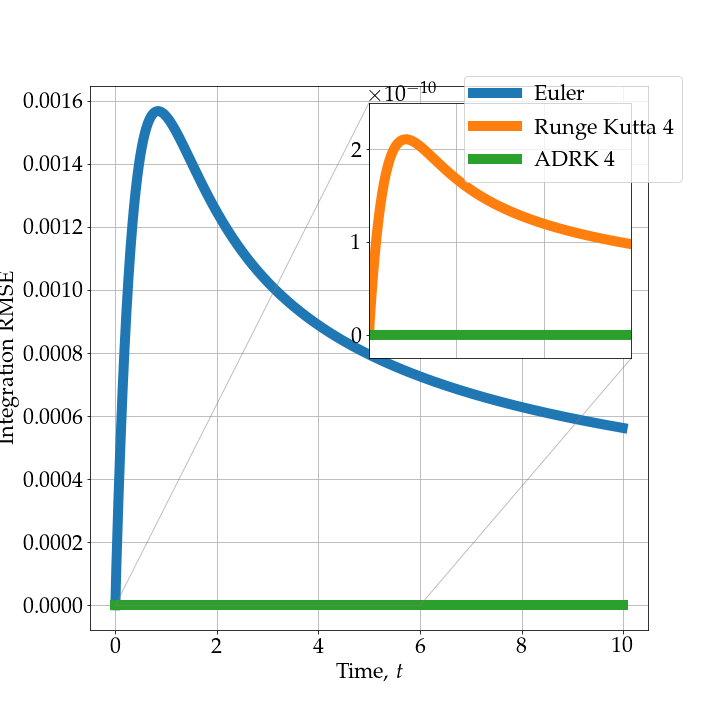}%
\caption{{{{ \bf  \em Root Mean Squared Error (RMSE) of several numerical resolutions of (\ref{eq:C3_exp_lin_ode}) with respect to the analytical solution}. The numerical resolution is carried given the same initial condition and integration time step for all the integration schemes.}
}}
\label{fig:c3_RMSE_pred_linear}
\end{figure}

We apply the training framework introduce in Section \ref{sec:ADRK_Non-linear_ODE_integration}, which relies on trajectory matching training loss. We report the forecasting performance of the trained model, including a comparison to the one of a Euler and RK4 scheme in figure \ref{fig:c3_RMSE_pred_linear}.

\begin{figure}[h]
\centering
  \includegraphics[clip,width=0.5\columnwidth,height=9cm]{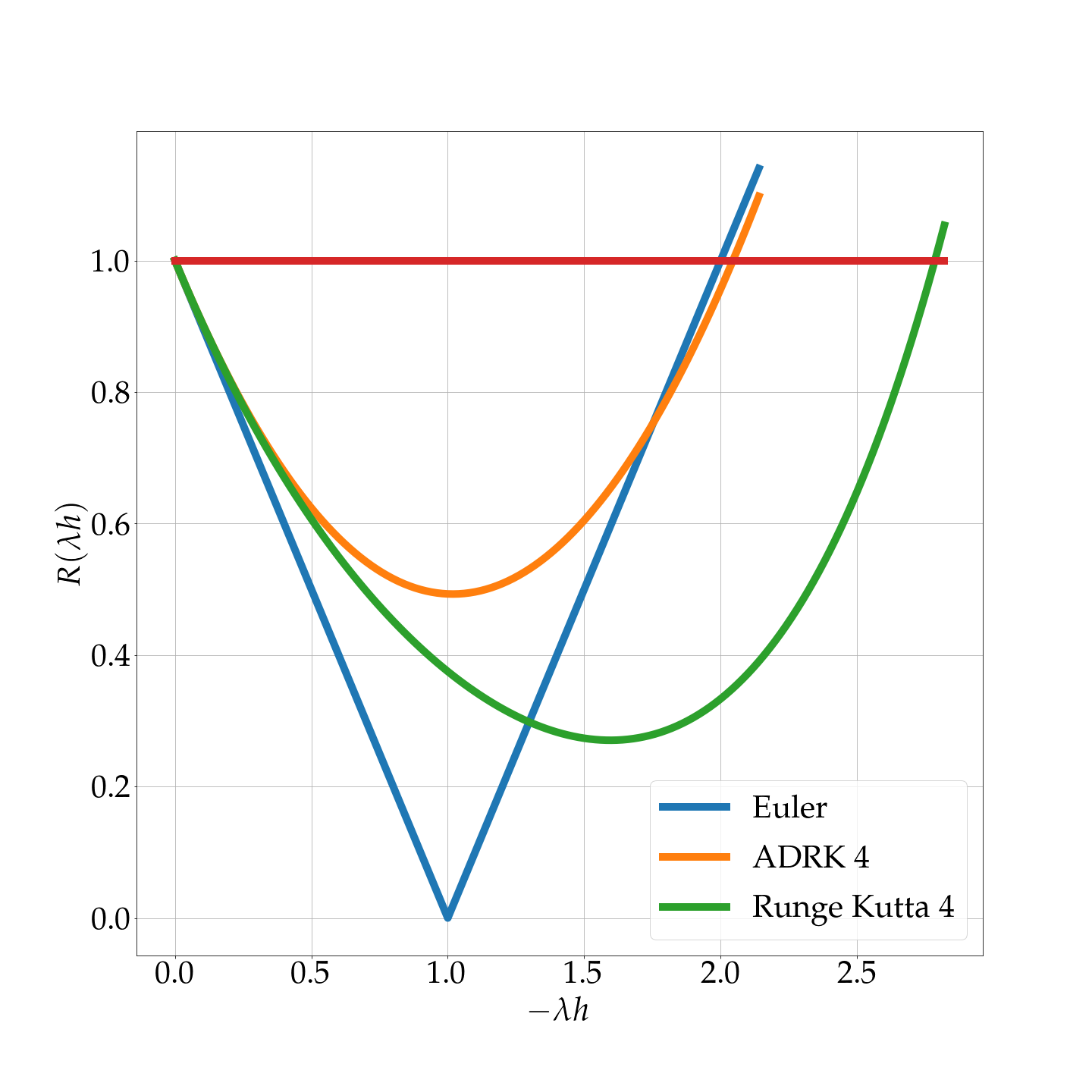}%
\caption{{{{ \bf  \em One-dimensional stability region of the proposed ADRK integration scheme with respect to classical state-of-the-art algorithms}. We compare the stability region of the proposed $\mathcal{ADRK}_4$ optimized on a short term forecasting cost with respect to classical state-of-the-art Euler and Runge-Kutta-4 integration schemes.}
}}
\label{fig:c3_stabil_region_1D}
\end{figure}

Besides, the forecasting performance, we analyse the stability region of the trained $\mathcal{ADRK}_4$ scheme as depicted in 
Fig. \ref{fig:c3_stabil_region_1D} with respect to state-of-the-art techniques.
The algorithms is stable for the integration time-step of the training data ($\lambda h = 0.2$). Our integration scheme also generalizes to a range of values of $\lambda h$ (up to $\lambda h = 2$) even if they were never used in the training phase. This is due to the integration scheme constraints (\ref{EQ:ADRK-q}) that were forced during the training phase.
Furthermore, we may also deduce the order of the trained $\mathcal{ADRK}_4$ from the analytical form of the stability polynomial $\mathcal{R}_{\mathcal{ADRK}_4}(\lambda h)$.
We can derive the stability polynomial of the trained ADRK as follows:
\begin{align}
\label{EQ:C3_Taylor_gain_RINN}
\begin{split}
\mathcal{R}_{\mathcal{ADRK}_4}(\lambda h)
              =& 1 + \lambda h + 0.4990 (\lambda h)^{2} +0.0037 (\lambda h)^{3} + 1.39\times 10^{-6}(\lambda h)^{4}
\end{split}
\end{align}
From the comparison of the polynomial (\ref{EQ:C3_Taylor_gain_RINN}) to the expansion of the true solution (\ref{EQ:C3_Taylor_Exp}), we determine that the trained ADRK scheme is a second-order integration scheme. 



\section{Non-linear ODE case-study}
\label{sec:Exp_L63}

In this section we report numerical experiments for the data-driven identification of Lorenz-63 dynamics. The Lorenz-63 dynamical system is a 3-dimensional model governed by the following ODE:
\begin{equation}
\label{eq:C3_lorenz_63}
\left \{\begin{array}{ccl}
\frac{d\vect{z}_{t,1}}{dt} &=&\sigma \left (\vect{z}_{t,2}-\vect{z}_{t,2} \right ) \\
\frac{d\vect{z}_{t,2}}{dt}&=&\rho \vect{z}_{t,1}-\vect{z}_{t,2}-\vect{z}_{t,1}\vect{z}_{t,3} \\
\frac{d\vect{z}_{t,3}}{dt} &=&\vect{z}_{t,1}\vect{z}_{t,2}-\beta \vect{z}_{t,3}
\end{array}\right.
\end{equation}
Under parameterization $\sigma =10$, $\rho=28$ and  $\beta=8/3$, this system involves chaotic dynamics with a strange attractor \cite{lorenz_deterministic_1963}.

We simulate Lorenz-63 state sequences using the LOSDA ODE solver \cite{odepack}. The integrated time-series was then sub-sampled based on a regular sampling rate. This state sequence is used as training data. In this section, we consider two experiments: i) the learning of an ADRK scheme to best reproduce Lorenz-63 trajectories knowing (\ref{eq:C3_lorenz_63}) and its parameterization
ii) the joint data-driven identification of  the underlying ODE and of an ADRK scheme Lorenz model to best reproduce the trajectories of the training dataset.

\subsection{ADRK scheme for Lorenz-63 system Integration}
\label{sec:ODE_integration_nlinear}
In this experiment, we compare a 4-stage ADRK scheme to the classical Runge-Kutta 4 method ($\mathcal{RK}_4$) for the integration of the Lorenz 63 system. The $\mathcal{ADRK}_4$ is optimized to match a single Lorenz 63 trajectory sampled at $h = 0.1, 0.15, 0.16, 0.17, 0.18, 0.19$. In table \ref{Tab:Integration_perfs_of_NL_ODE}, we report the results of this experiment by specifying the ability of the $\mathcal{RK}_4$ and the $\mathcal{ADRK}_4$ schemes to correctly integrate the Lorenz 63 equation at these integration time-steps. Specifically, an integration scheme succeeds if it generates chaotic Lorenz 63 trajectories and fails if it makes the integration blows up, generates stable equilibrium points or other non-chaotic limit-sets.

\begin{table}
\centering
\begin{tabular}{@{} l *6c @{}}
\toprule
 \multicolumn{1}{c}{h}    & 0.1  & 0.15  & 0.16  & 0.17 & 0.18& 0.19  \\ 
\midrule
\midrule
 $\mathcal{ADRK}_4$ &  \checkmark & \checkmark & \checkmark & \checkmark & \text{\sffamily X} & \text{\sffamily X}\\ 
 $\mathcal{RK}_4$   & \checkmark & \text{\sffamily X} & \text{\sffamily X} & \text{\sffamily X}&\text{\sffamily X}&\text{\sffamily X} \\
 \bottomrule
 \end{tabular}
\caption{{{{ \bf  \em Integration ability of both the $\mathcal{RK}_4$ and the trained $\mathcal{ADRK}_4$ on the Lorenz 63 system with different integration time steps}.}
}}
 \label{Tab:Integration_perfs_of_NL_ODE}
\end{table}

\begin{figure*}
\centering
  \subfloat[]{%
\includegraphics[width=0.3\textwidth]{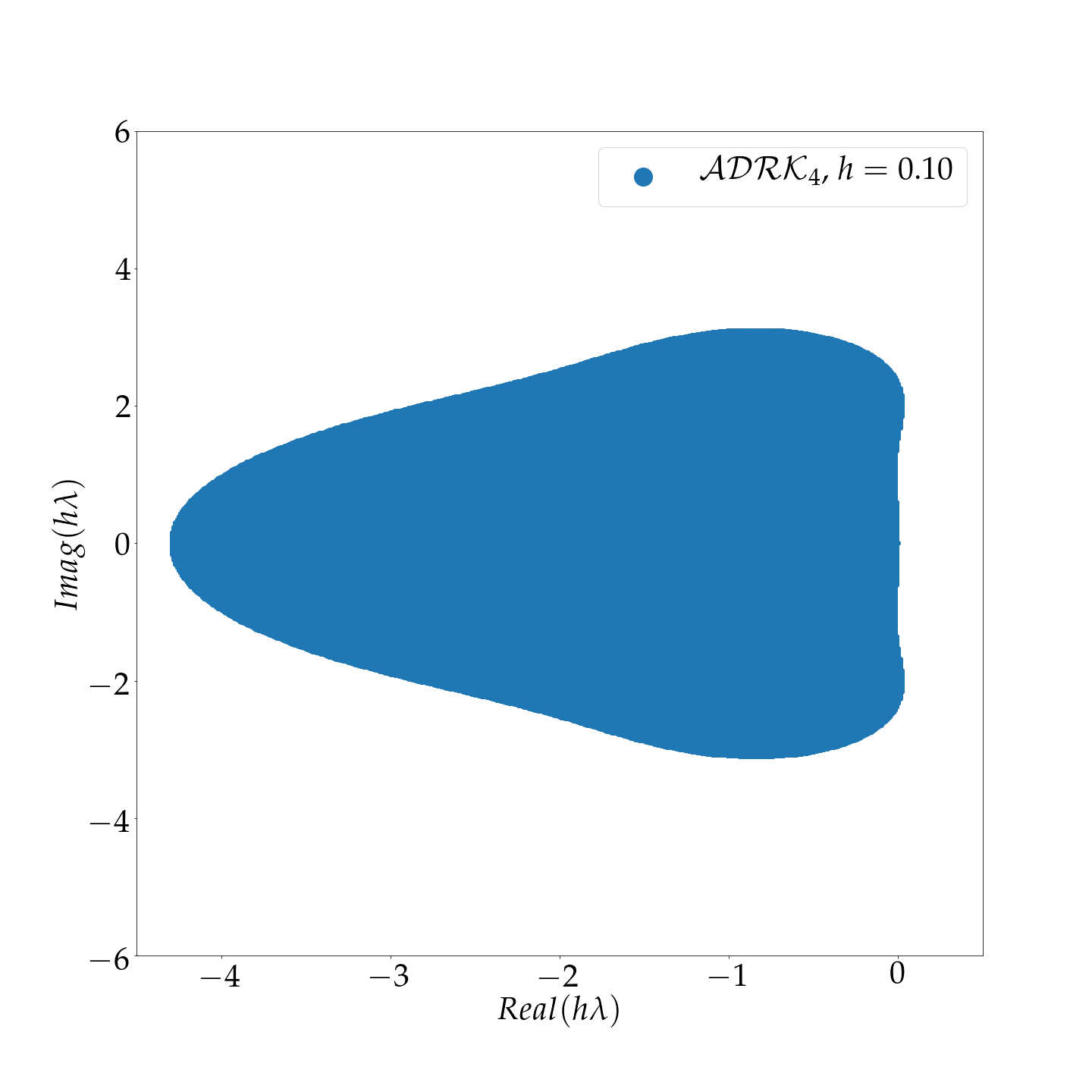}}
  \subfloat[]{%
\includegraphics[width=0.3\textwidth]{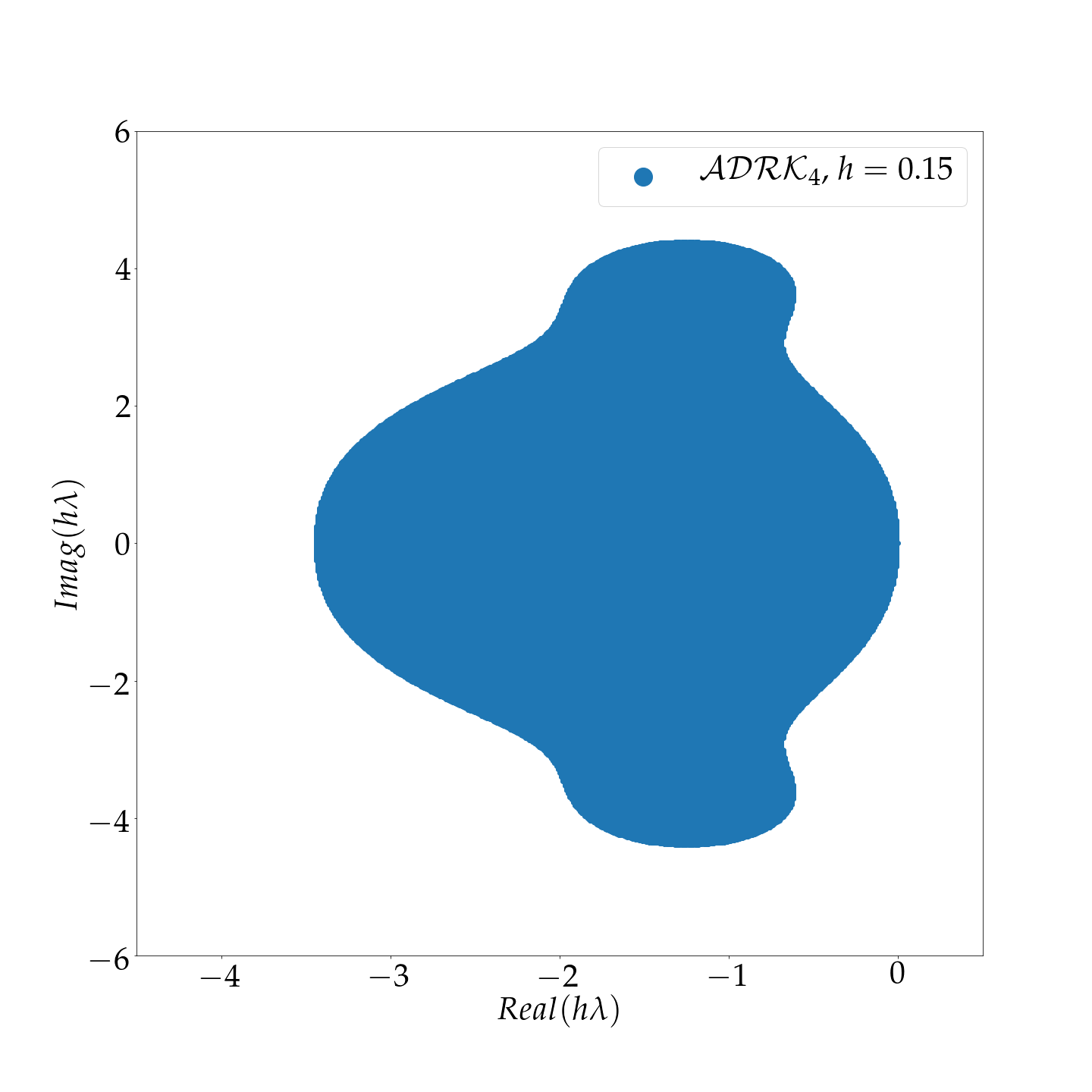}}
  \subfloat[]{%
\includegraphics[width=0.3\textwidth]{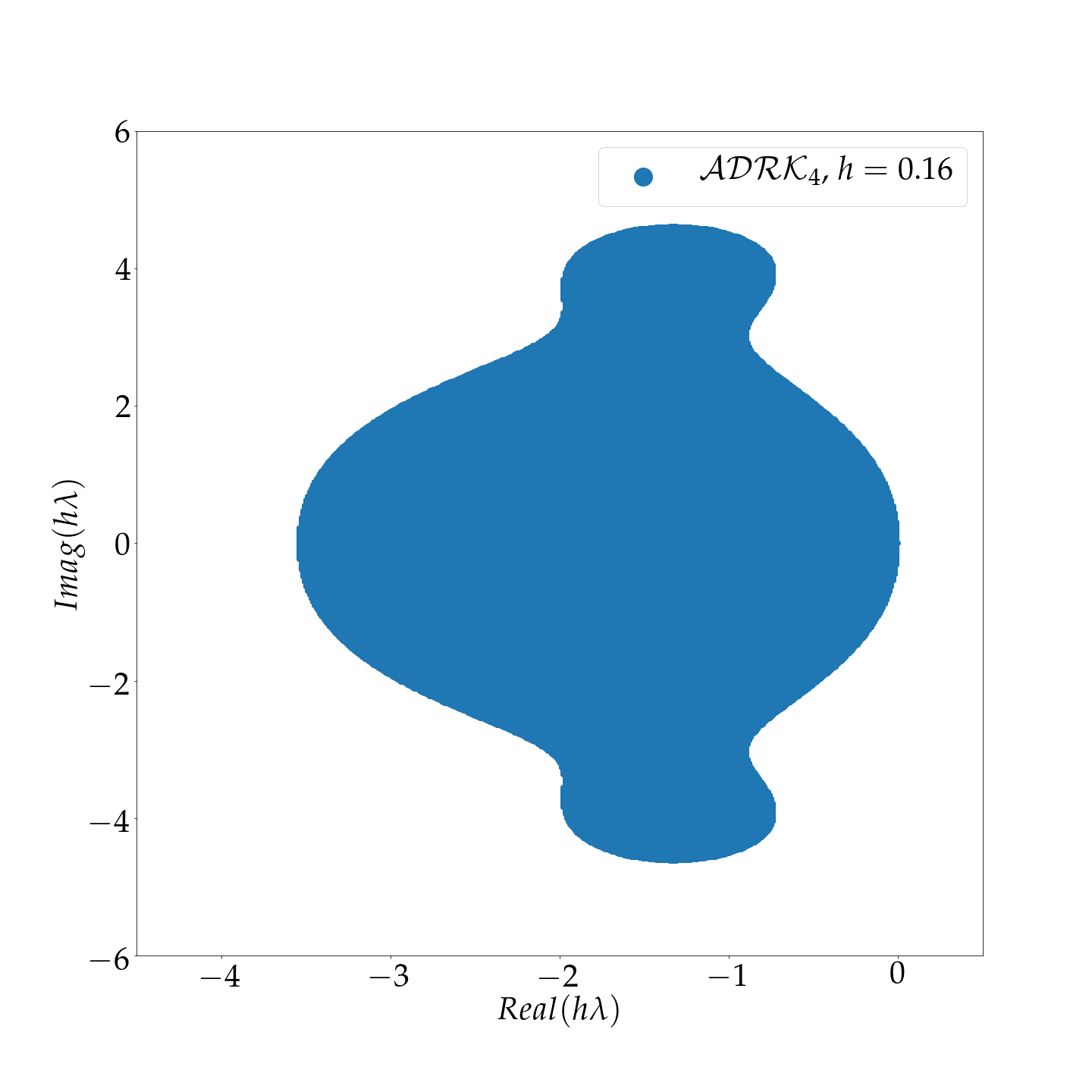}}\\
  \subfloat[]{%
\includegraphics[width=0.3\textwidth]{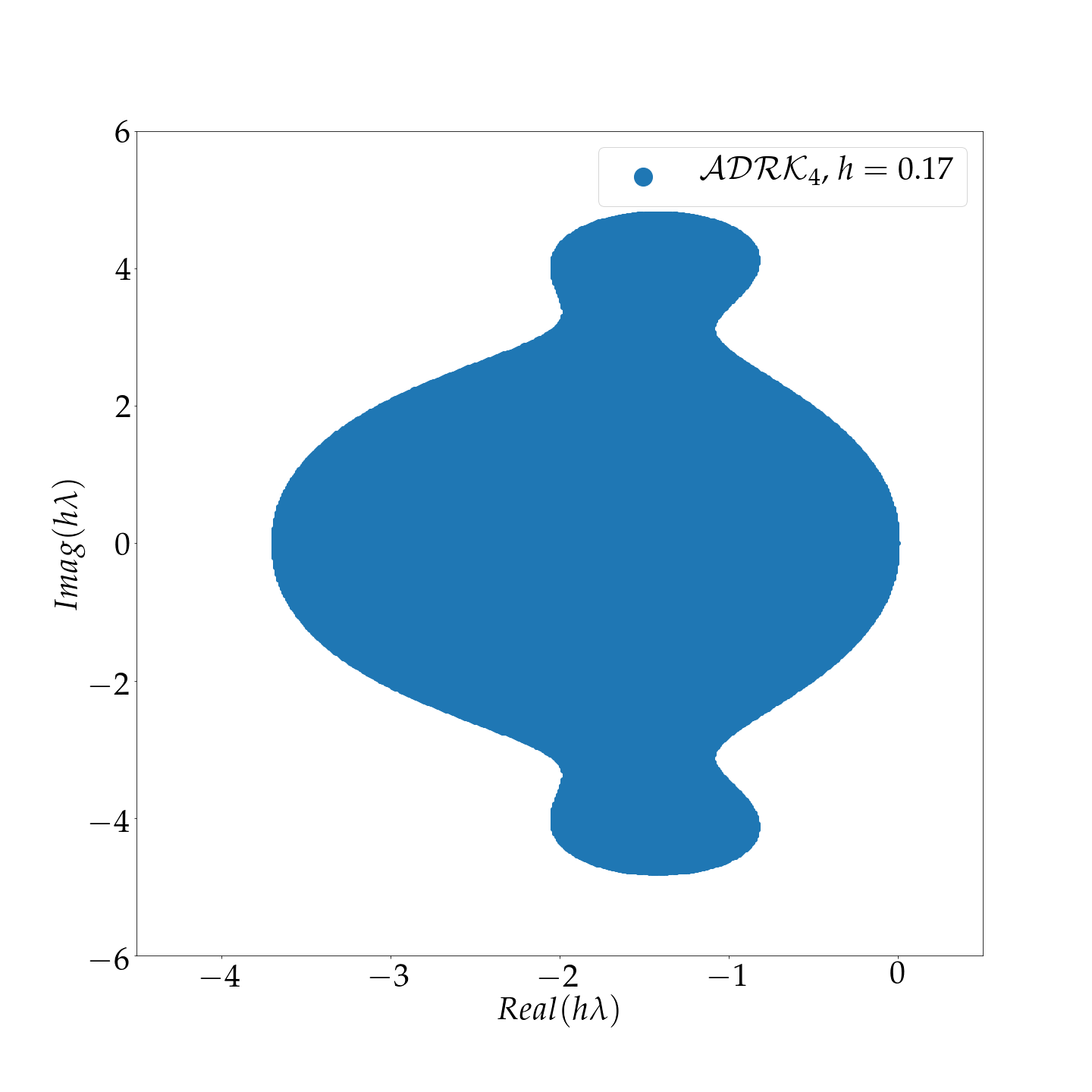}}
  \subfloat[]{%
\includegraphics[width=0.3\textwidth]{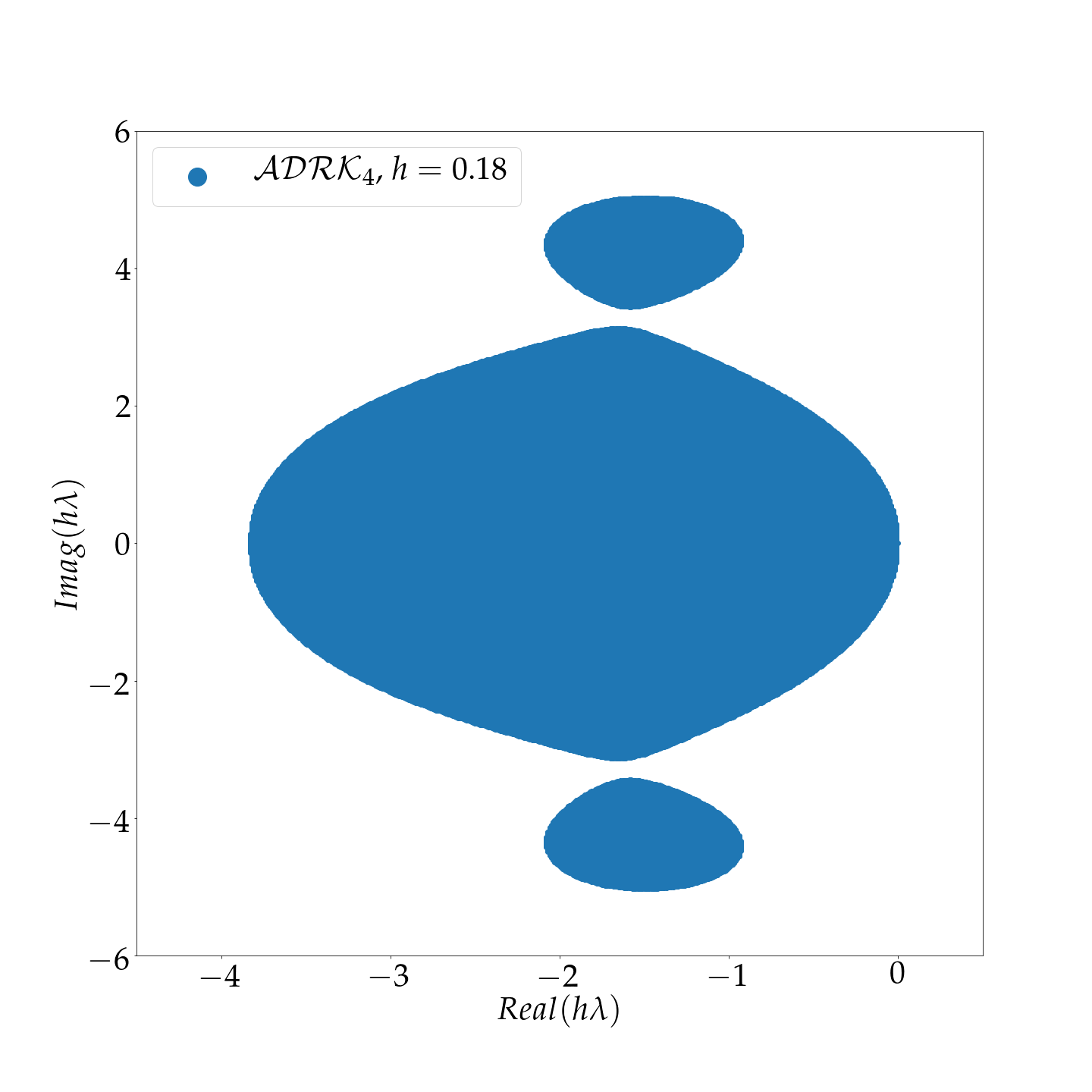}}
  \subfloat[]{%
\includegraphics[width=0.3\textwidth]{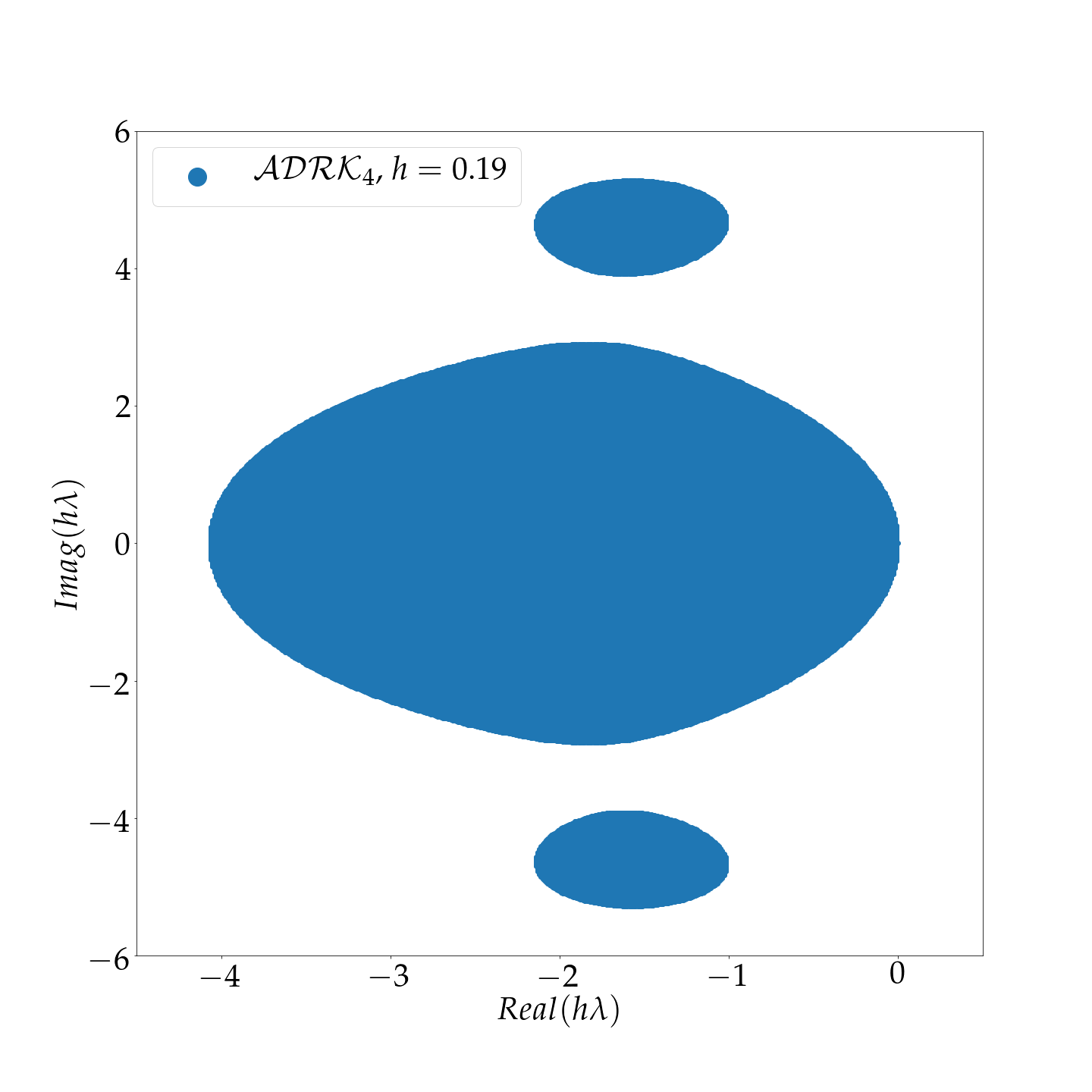}}
\\
  \subfloat[]{%
\includegraphics[width=0.3\textwidth]{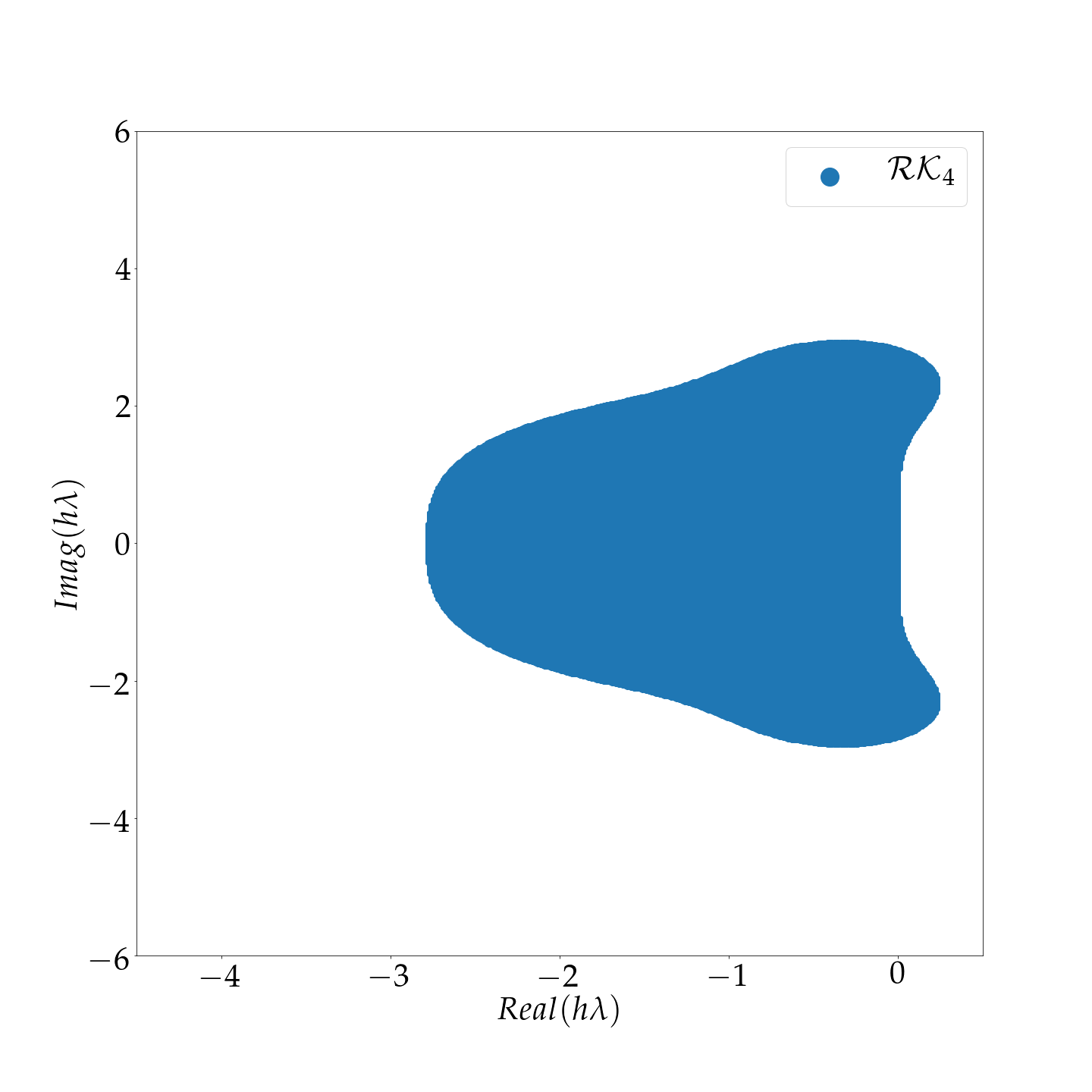}}
\caption{{{{ \bf  \em Stability region of the $\mathcal{ADRK}_4$ schemes trained to match the Lorenz 63 states at different integration time-steps}. Stability region of the $\mathcal{ADRK}_4$ scheme trained based on an integration time-step ranging from $h = 0.1$ to $h = 0.19$ are given in figures (a) to (f) respectively. The stability region of the classical $\mathcal{RK}_4$ is given as a reference in fig. (g).}
}}
\label{fig:Stability_polynom_integrationL63}
\end{figure*}

The $\mathcal{ADRK}_4$ is able to correctly generate chaotic trajectories from the Lorenz 63 ODE up to a sampling rate $h = 0.17$. By contrast, the classical Runge-Kutta 4 technique either blows up or get stuck at an equilibrium point for $h \ge 0.15$. 
Whereas the latter is derived independently of any non-linear ODE to match a fourth order truncation of the Taylor expansion (\ref{EQ:C3_Taylor_Exp}) leading to a smaller acceptable integration time-step for the Lorenz equation, 
the good properties of the trained $\mathcal{ADRK}_4$ scheme comes from the data-driven training procedure. Figure \ref{fig:Stability_polynom_integrationL63} illustrates the stability region of the trained ADRK schemes on each integration time-step. When compared to the classical RK(4) scheme, the ADRK schemes involve larger stability regions as the integration time-step $h$ increases from $0.1$ to $0.17$, which makes them relevant for the integration of Lorenz 63 ODE. 

Interestingly, the stability region mostly increase in the imaginary axis for the Lorenz equation which emphasises that this equation is probably dominated by an advection process. This observation can be explained through a linearization of the Lorenz 63 model based on the Koopman operator \cite{koopman1931hamiltonian}. Although such a linearization would theoretically lead to an infinite dimensional system, an approximate finite dimensional model can be written as proposed in \cite{brunton2017chaos} with an additional forcing term. In such representation, the linear dynamical system, shadowing the true non-linear Lorenz-63 equation, will have only imaginary eigenvalues. This property can explain the stability growth towards the imaginary axis rather than the real one. The $\mathcal{ADRK}_4$ scheme is able to match these requirements leading to a larger acceptable step-size for Lorenz-63 system. 

Overall, this experiment highlights one of the main potentials of the proposed framework when considering the time integration of non-linear differential equations. In these applications, and especially when considering complex problems for which a detailed analytical study of the underlying dynamical properties can hardly be achieved, the proposed ADRK framework adapts to the observed non-linear dynamics leading to an equation-adapted integration scheme.

\subsection{Lorenz-63 system Identification}
\label{sec:identification_L63}
In this experiment we subsampled the Lorenz 63 trajectories to a low regular sampling rate $h_1 = 0.2$, $h_2 = 0.3$ and $h_3 = 0.4$. The goal of this experiment is to try to discover a model for the Lorenz system given temporally sparse data that can not be integrated using classical integration schemes. For this purpose, an approximate model $f_{\theta_{NN}}$ given in (\ref{EQ:sec3_ODE_eq_DD}) is optimized to generate a Lorenz state sequence sampled at a sparse sampling rate. The approximate model is a Linear Quadratic Model (LQM) as proposed in \cite{fablet_bilinear_2017}. This architecture ensures that the true model lies within the space of possible model parameterizations. For benchmarking purpose, the following models were tested:
\begin{itemize}[topsep=0pt]
\itemsep0em 
\item \textbf{Sparse regression model} \cite{brunton_discovering_2016} (SR): This model computes a sparse regression over an augmented states vector based on second order polynomial representations of the Lorenz states. The learnt dynamical model is then integrated to compute forecasts using the LOSDA ODE solver \cite{odepack}. 
\item \textbf{Euler based model} (Euler): The approximate model $f_{\theta_{NN}}$ is integrated using an explicit Euler scheme.

\item \textbf{Runge-Kutta 4 based model} ($\mathcal{RK}_4$): The approximate model $f_{\theta_{NN}}$ is integrated using the classical Runge-Kutta 4 scheme. 

\item \textbf{Neural ODE} (DOPRI8): The approximate model $f_{\theta_{NN}}$ is integrated using the adaptive step size dopri8 solver. The backward pass is computed here using the adjoint sensitivity method \cite{alma990001535520306161}, as proposed in \cite{chen2018neural}. 

\item \textbf{11-stage ADRK scheme}: Proposed automatic differentiation Runge Kutta scheme with a number of stages equal to 11. In this architecture, the weights of the integration scheme are learnt jointly with the parameters of the approximate model as explained in section \ref{sec:ADRK_Identification}.

\end{itemize}
Besides the proposed ADRK scheme, all the benshmark models use a specified integration routines and the optimization is carried only with respect to the parameters of the model $f_{\theta_{NN}}$. Finally the data-driven integration schemes are noted $\mathcal{ADRK}_{h_1}$, $\mathcal{ADRK}_{h_2}$ and $\mathcal{ADRK}_{h_3}$ with each index ($h_1$, $h_2$ and $h_3$) corresponding to the time sampling of the training time series.

\begin{table}
\centering
\begin{tabular}{ll*{5}c}
\toprule
\multicolumn{2}{c}{Model} & $h=0.2$ & $h=0.3$ & $h=0.4$ \\
\midrule \midrule 
\multirow{2}{*}{SR}
& $t_0+h$  & $9.06$ & $>10$  & $9.34$ \\
& $t_0+4h$ & $6.02$ & $5.81$ & $6.97$  \\ 
\midrule
\multirow{2}{*}{Euler}
& $t_0+h$  & $4.27$ & $2.57$  & $1.99$ \\
& $t_0+4h$ & $>10$ & $>10$ & $7.89$    \\ 
\midrule
\multirow{2}{*}{$\mathcal{RK}_4$}
& $t_0+h$  & $2.05$ & $3.10$  & $2.48$  \\
& $t_0+4h$ & $3.82$ & $7.33$  & $>10$    \\ 
\midrule 
\multirow{2}{*}{Dopri8}
& $t_0+h$  &  \textbf{0.005} & \textbf{0.0001} & 3.1305 \\
& $t_0+4h$ & {0.021}         & \textbf{0.0003} & $>10$  \\ 
\midrule 
\multirow{2}{*}{$\mathcal{ADRK}$}
& $t_0+h$ &  ${0.017}$ & ${0.077}$ & $\textbf{0.189}$ \\
& $t_0+4h$ & $\textbf{{0.020}}$ & ${0.23}$ & $\textbf{1.93}$ \\ 
\bottomrule
\end{tabular}
\caption {{{ \bf  \em Forecasting performance of data-driven models for Lorenz-63 dynamical model.} Mean Root Mean Squared Error (RMSE) for different forecasting time-steps of the tested models.}}
\label{tab:C3_forecast_Lorenz63}
\end{table}

\begin{figure}[h]
\centering
  \includegraphics[clip,width=1\columnwidth,height=12cm]{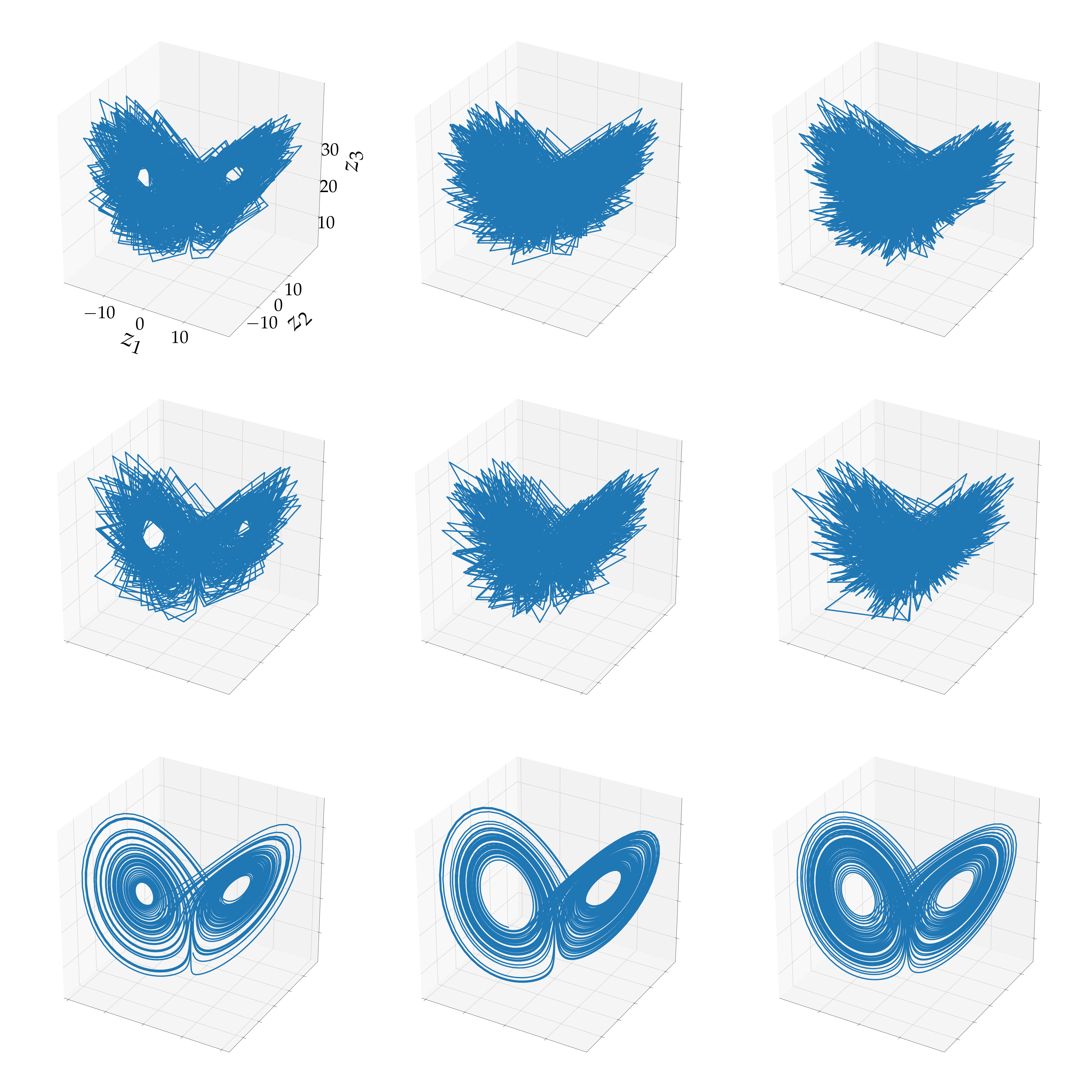}%
\caption{{{{ \bf  \em Training data and simulated attractors from the corresponding data-driven $\mathcal{ADRK}$ based models}. First row, training sequence generated with different regular time sampling; second row, simulated attractors of the data-driven models, third row; simulated attractors using a smaller integration time-step. The columns correspond to the sampling rate of the training series (ranging from 0.2 to 0.4.) which corresponds to the integration time-step $h$ when considering fixed step-size algorithms.}
}}
\label{fig:c3_Attractors_several_h}
\end{figure}

\begin{figure}[h]
\centering
  \includegraphics[clip,width=1\columnwidth,height=12cm]{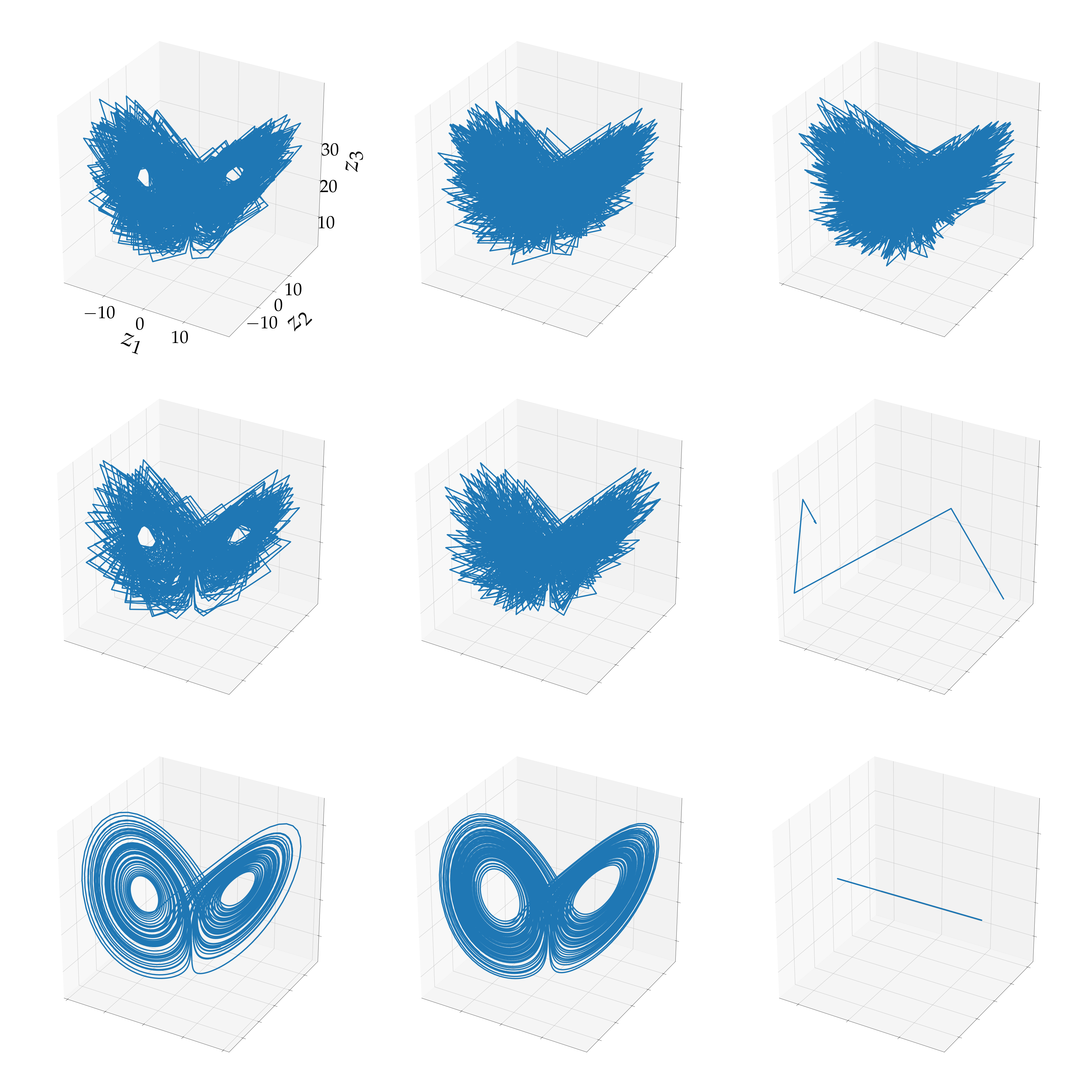}%
\caption{{{{ \bf  \em Training data and simulated attractors from the corresponding data-driven Dopri8 based models}. First row, training sequence generated with different regular time sampling; second row, simulated attractors of the data-driven models, third row; simulated attractors using a smaller integration time-step. The columns correspond to the sampling rate of the training series (ranging from 0.2 to 0.4.).}
}}
\label{fig:c3_Attractors_several_h_dopri}
\end{figure}

\textbf{ADRK Performance analysis}: We report the forecasting performance in Table \ref{tab:C3_forecast_Lorenz63}. Figures \ref{fig:c3_Attractors_several_h} and \ref{fig:c3_Attractors_several_h_dopri} illustrate the trajectories generated using the trained $\mathcal{ADRK}$ and Dopri8 based models respectively. When compared to classical fixed stipe-size models, the proposed $\mathcal{ADRK}$ is the only model able to discover the hidden dynamics of the system. The reason is that classical fixed step-size solvers such as Euler and Runge-Kutta 4 lead to high truncation errors making the identifiability impossible, the proposed model in the other hand, when deployed with 11 stages, can mimic adapted, stable and precise integration schemes, that are able to unfold the true Lorenz 63 dynamics. These results can be explained by analysing the truncation error of a q-stage ADRK scheme. The truncation error of numerical integration scheme provides a straightforward tool to state about the one step ahead integration performance of a given integration technique as a function of its order $p$. Assuming that the $q$-stage ADRK corresponds to a $\hat{p}$-order numerical integration scheme, the loss function of the ADRK relates, following (\ref{EQ:sec3_Truncation_error}), to the truncation error of the learnt integration scheme as follows:
\begin{equation}
\label{EQ:sec3_Truncation_error_training}
\mathcal{L}(\Theta) \propto (\epsilon_n)^2 = (E(\hat{p},f_{\theta_{NN}},\vect{z}_{t_{n}},t_{n})h^{p+1})^2
\end{equation}
where $E \in \mathbb{R}$ a function of $\hat{p},f_{\theta_{NN}},\vect{z}$ and $t_{n}$.



Equation (\ref{EQ:sec3_Truncation_error_training}) illustrates a main characteristic about learning data-driven representations of dynamical models, using fixed step-size integration schemes: (ii) we may jointly tune $f_{\theta_{NN}}$ and $\hat{p}$ in the ADRK architecture to lower the training loss function.

Assuming that the integration time-step $h$ is set by the temporal sampling of our training data, one may decrease the training error through tuning the approximation $f_{NN}$ of the true dynamical model $f$\cite{brunton_discovering_2016,fablet_bilinear_2017}. One may also decrease the training loss function through an adapted order $\hat{p}$ of the integration scheme, parameterized in the proposed ADRK framework by the number of stages $q$. Overall, the training loss lower bound motivates the use of the proposed framework in identification applications since if we start an identification problem with a fixed integration scheme or family that can not integrate properly the dynamics (which is a fair enough claim since we are interested in non-linear system identification where stability and precision criteria are not known a priori), we might not find a good fit for our data simply due to the choice of the integration method.

\begin{figure*}
\centering
  \subfloat[]{%
\includegraphics[width=0.5\textwidth]{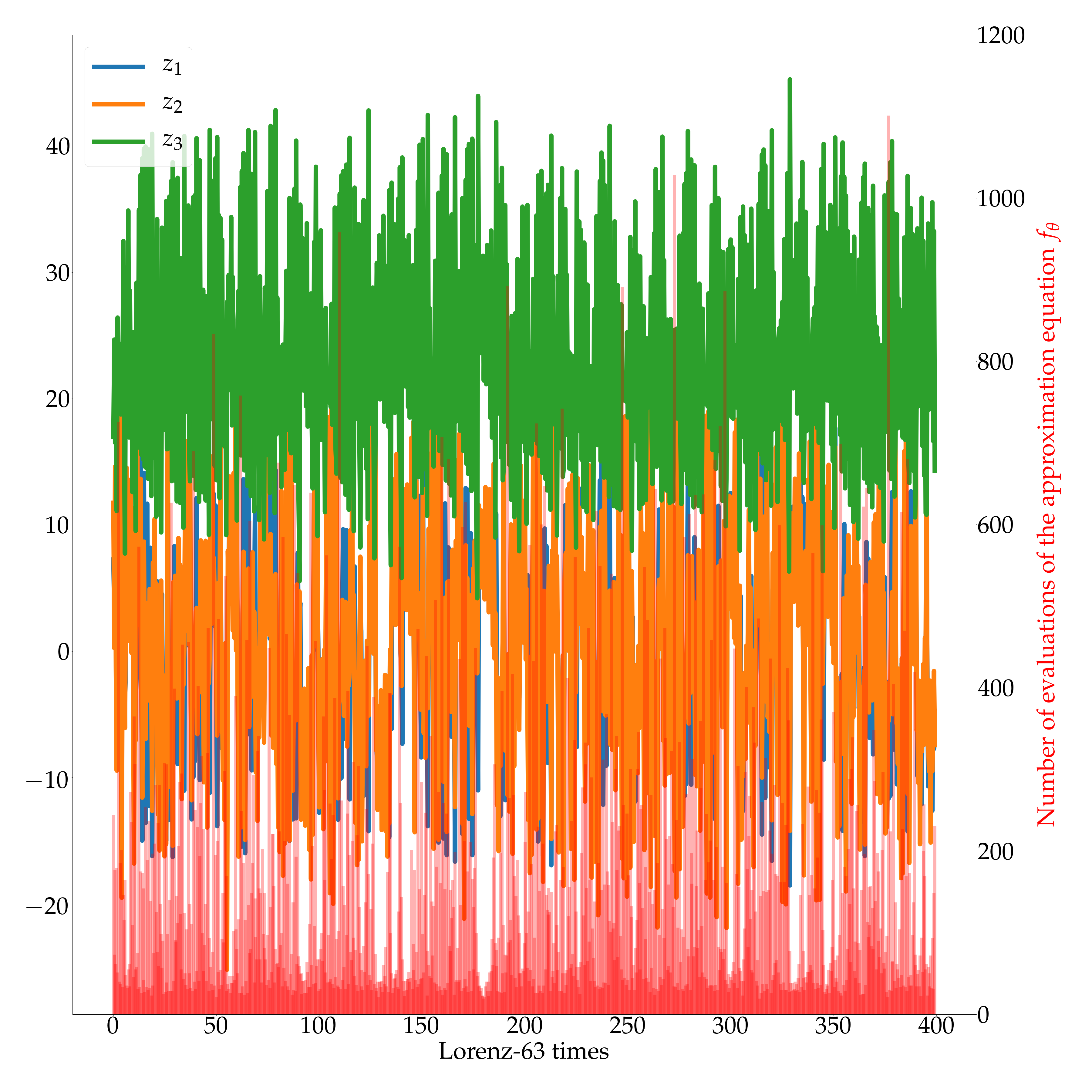}}
  \subfloat[]{%
\includegraphics[width=0.5\textwidth]{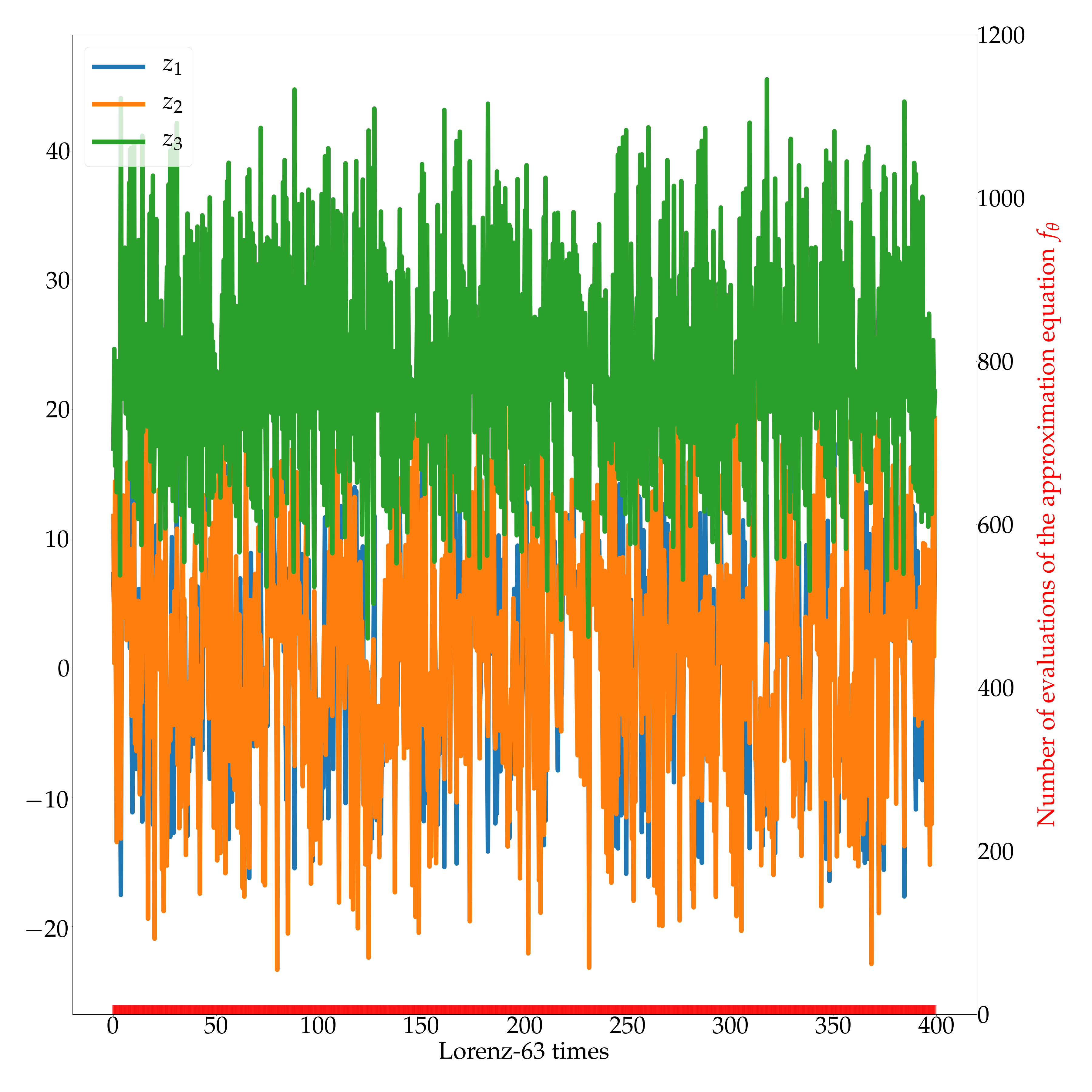}}
\caption{{{{ \bf  \em Simulated time series of the Lorenz 63 dynamics from two data-driven models}. (a) Dopri8 based model trained without the adjoint technique proposed in \cite{chen2018neural}. (b) Proposed $\mathcal{ADRK}$. The training of both these models was carried on the dataset sampled with $h = 0.4$. The red bars correspond, in the $\mathcal{ADRK}$ figure (b) to the exact number of evaluations of the approximate function $f_\theta$ for every point within the simulated trajectory. In the DOPRI8 figure, the bars correspond to the minimum number of evaluations, computed when considering each integration carried within a single stage integration scheme.}
}}
\label{fig:C3_nb_of_func_evaluation}
\end{figure*}

The poor results of the sparse regression method is in the other hand simply due to a wrong estimation of the derivatives. This step is inevitable using such technique and when provided with temporally sparse data, a decent estimation becomes impossible. Adaptive step-size based models in the other hand lead to an overall better short term forecast and, similarly to the proposed framework, are able to correctly identify the Lorenz 63 model when provided with data sampled at $h = 0.2$ and $h = 0.3$. This method however is unable to derive a decent approximation when provided with data sampled at $h = 0.4$. We believe that this is principally due to the backward pass formulation proposed in \cite{chen2018neural} which integrates an augmented ODE backward in time in the training phase (to avoid heavy memory usage due to classical backpropagation through residual network). These assumptions were then verified on the same case study by computing the backward pass through storing every evaluation of the adaptive solver. Within this configuration, the Dopri8 technique is able to unfold the true structure of the attractor and achieves a decent short term forecast, similar to the proposed $\mathcal{ADRK}$ (0.18 at $t_0+h$ and 0.94 at $t_0+4h$). Figure \ref{fig:C3_nb_of_func_evaluation} illustrates a simulation example carried using both the $\mathcal{ADRK}$ and the Dopri8 based models on the dataset sampled at $h = 0.4$ (the Dopri8 based model is trained without the use of the adjoint method). The number of function evaluations of the adaptive solver varies between 19 and 14313 (please refer to figure \ref{fig:C3_nb_of_func_evaluation} for a visualisation of a lower bound of the number of function evaluations per initial condition) against only 11 for the proposed framework making our method highly computationally efficient, especially when considering more complex systems. Furthermore, storing all the integration steps of an adaptive solver in order to backpropagate the training error is not guaranteed to be feasible in first place on more complex, high dimensional systems due to memory blowups.

\begin{figure*}
\centering
  \subfloat[]{%
\includegraphics[width=0.3\textwidth]{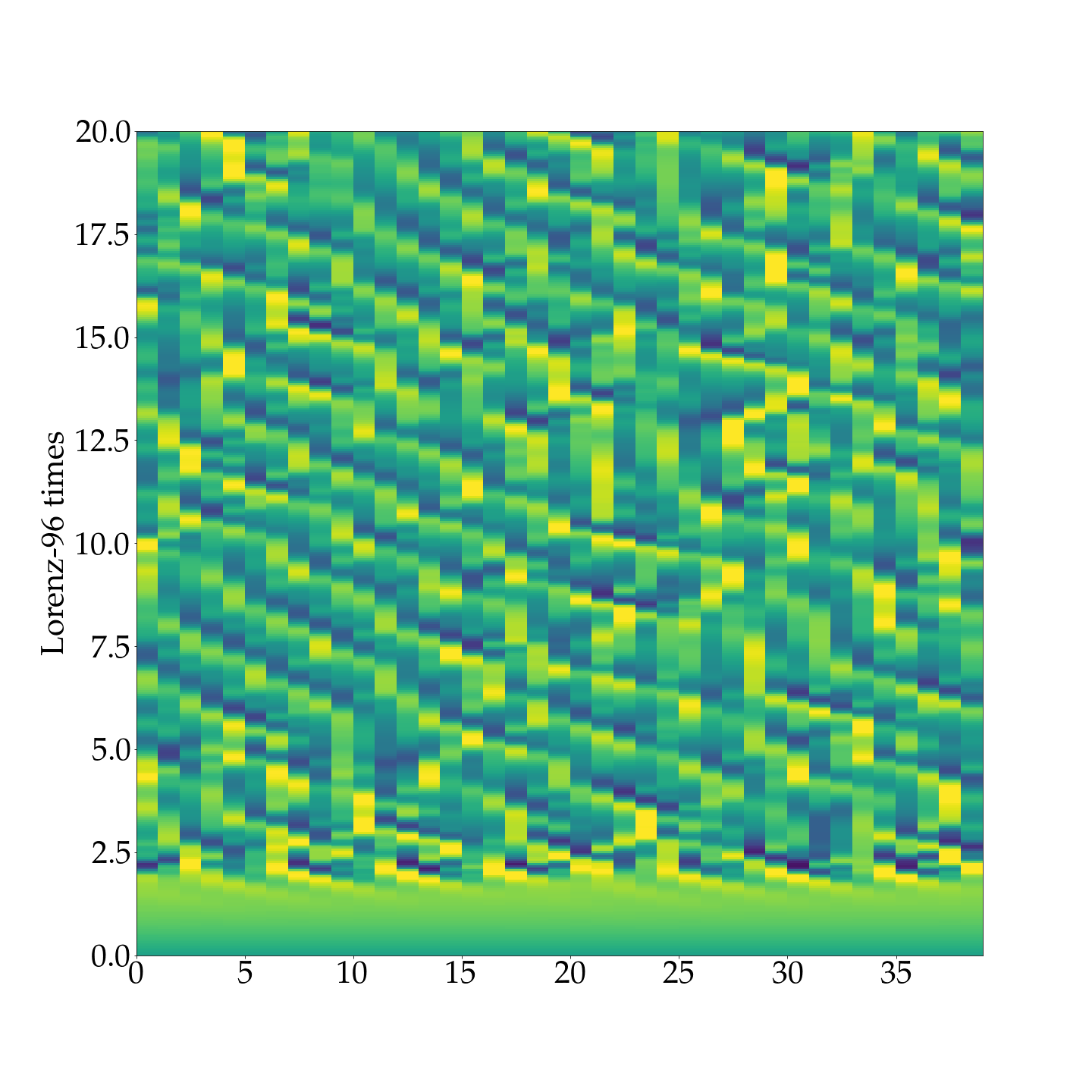}}
  \subfloat[]{%
\includegraphics[width=0.3\textwidth]{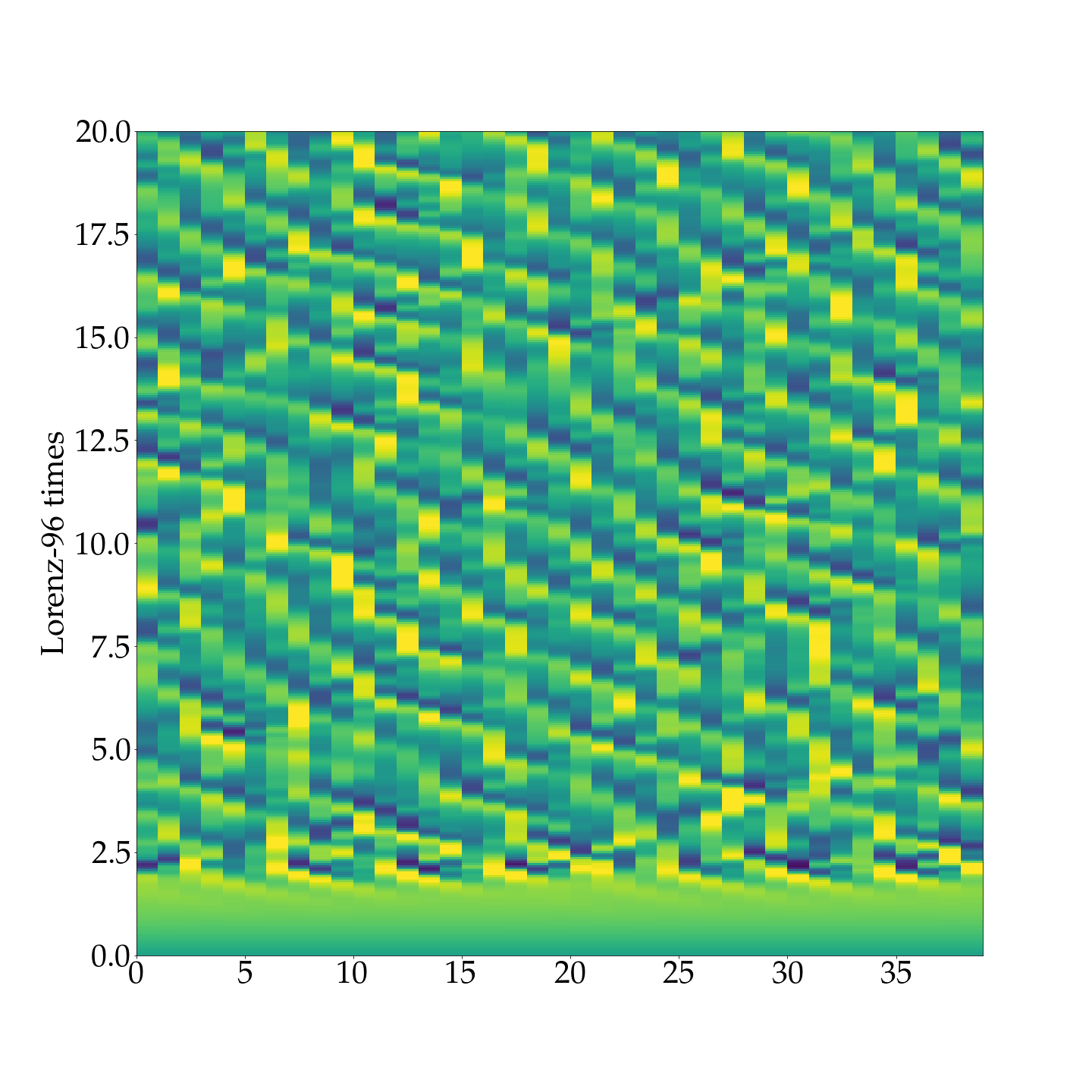}}
  \subfloat[]{%
\includegraphics[width=0.3\textwidth]{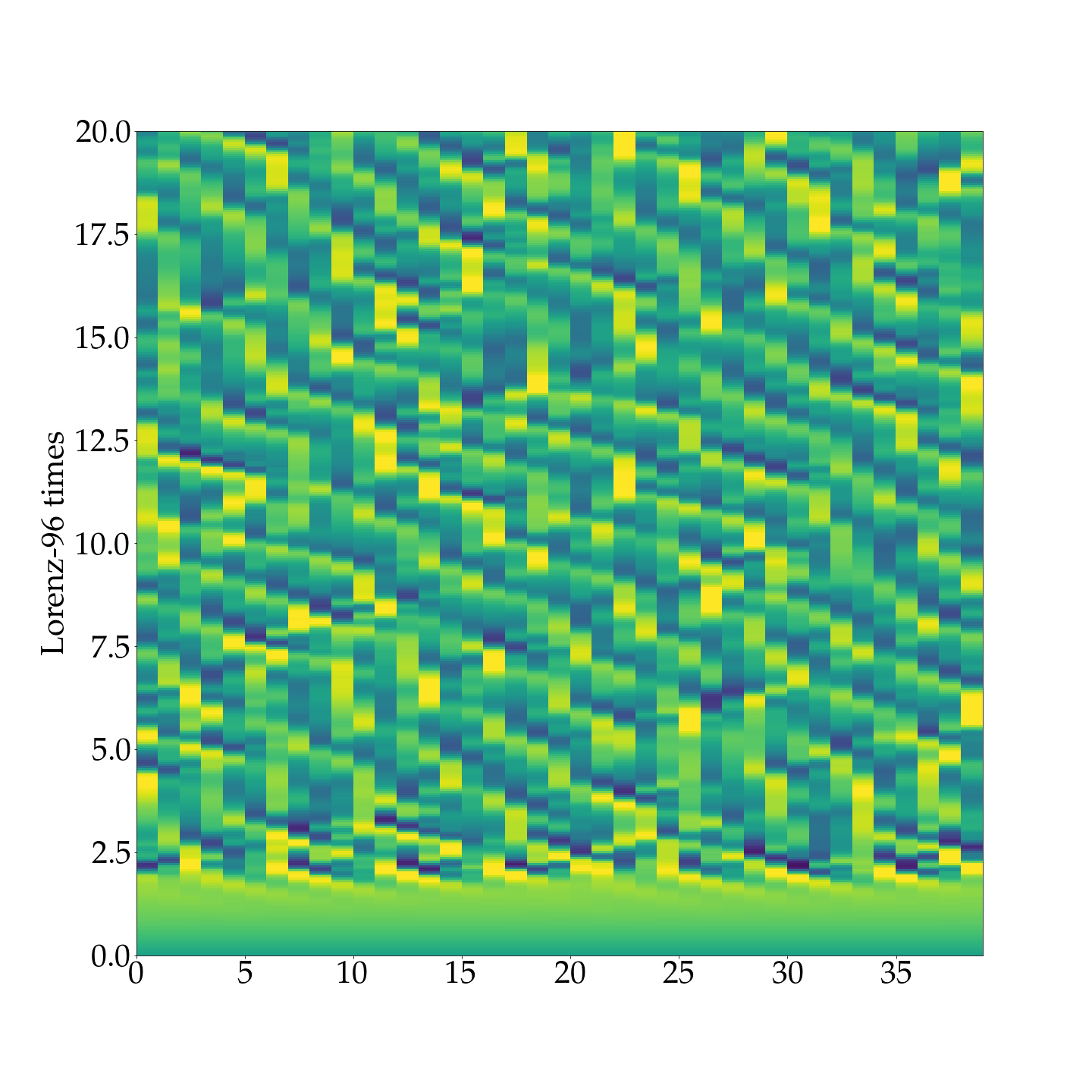}}
\caption{{{{ \bf  \em Application of the learnt integration scheme on the Lorenz 96 equation}. Integration of the Lorenz 96 model using the $\mathcal{ADRK}_{h_1}$ in (a), the $\mathcal{ADRK}_{h_2}$ in (b) and the $\mathcal{ADRK}_{h_3}$ in (c).}
}}
\label{fig:C3_RINN_L96_several_h}
\end{figure*}

\textbf{Generalizability of the learnt ADRK schemes to new integration configurations}: Another interesting experiment is illustrated in Fig. \ref{fig:C3_RINN_L96_several_h}. The trained integration scheme is dissociated from the data-driven model $f_{NN}$ and used to integrate the Lorenz 96 differential equation. Since the integration scheme is always consistent, there is a range of integration time-steps for which the integration scheme is stable given a new ODE and thus, converges to give a correct simulation of the equation.

\begin{figure}[h]
\centering
  \includegraphics[clip,width=0.5\columnwidth,height=9cm]{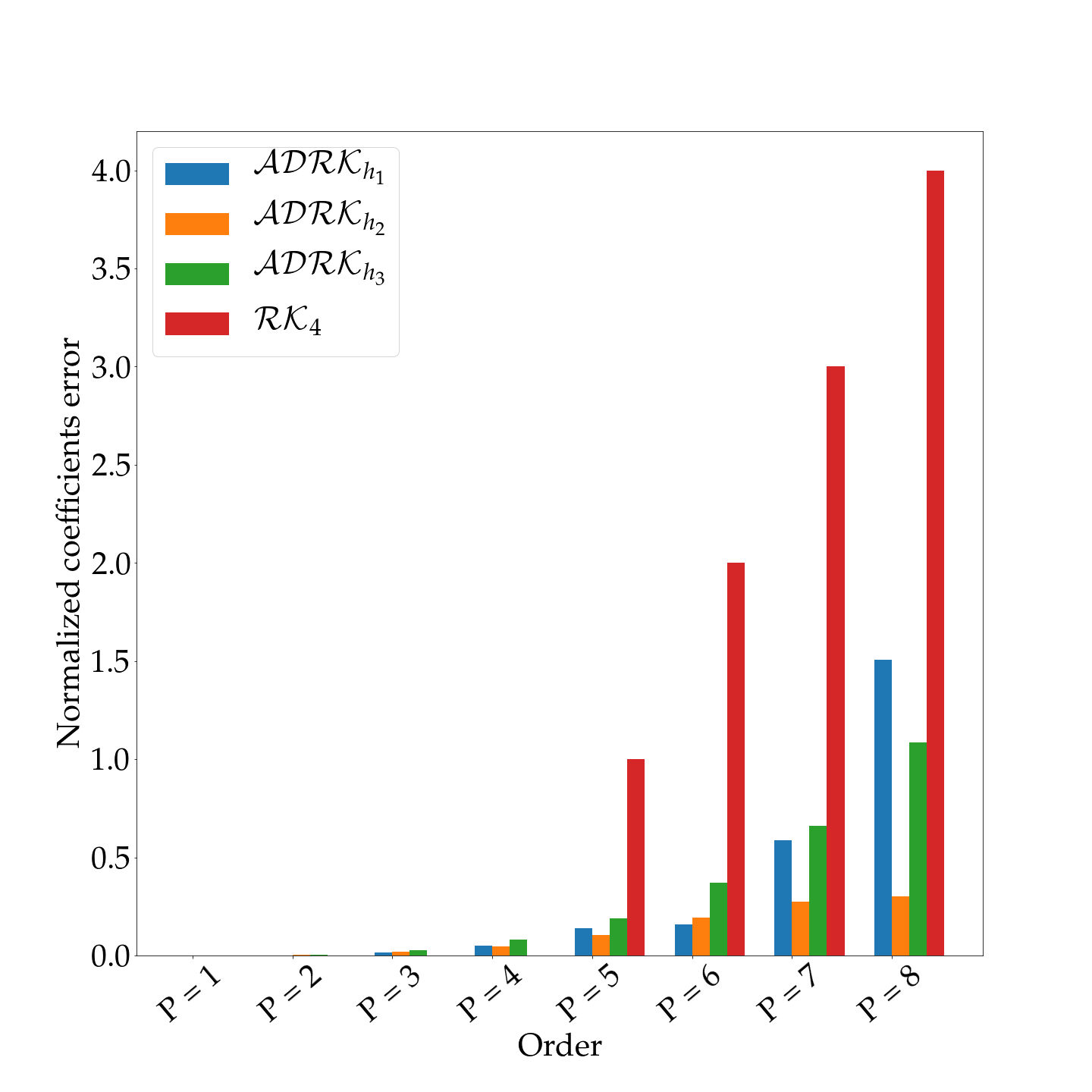}%
\caption{{{{ \bf  \em Normalized coefficients error of the stability polynomials of the ADRK schemes (Lorenz 63 identification case-study) with respect to the Taylor expansion of the analytical solution}. Normalized cumulative error of the $\mathcal{ADRK}$ stability polynomials coefficients up to the order p, with respect to the Taylor expansion of the true solution. The Runge-Kutta-4 algorithm is given as a reference.}
}}
\label{fig:c3_coeff_errors}
\end{figure}

\begin{figure}[h]
\centering
  \includegraphics[clip,width=0.5\columnwidth,height=9cm]{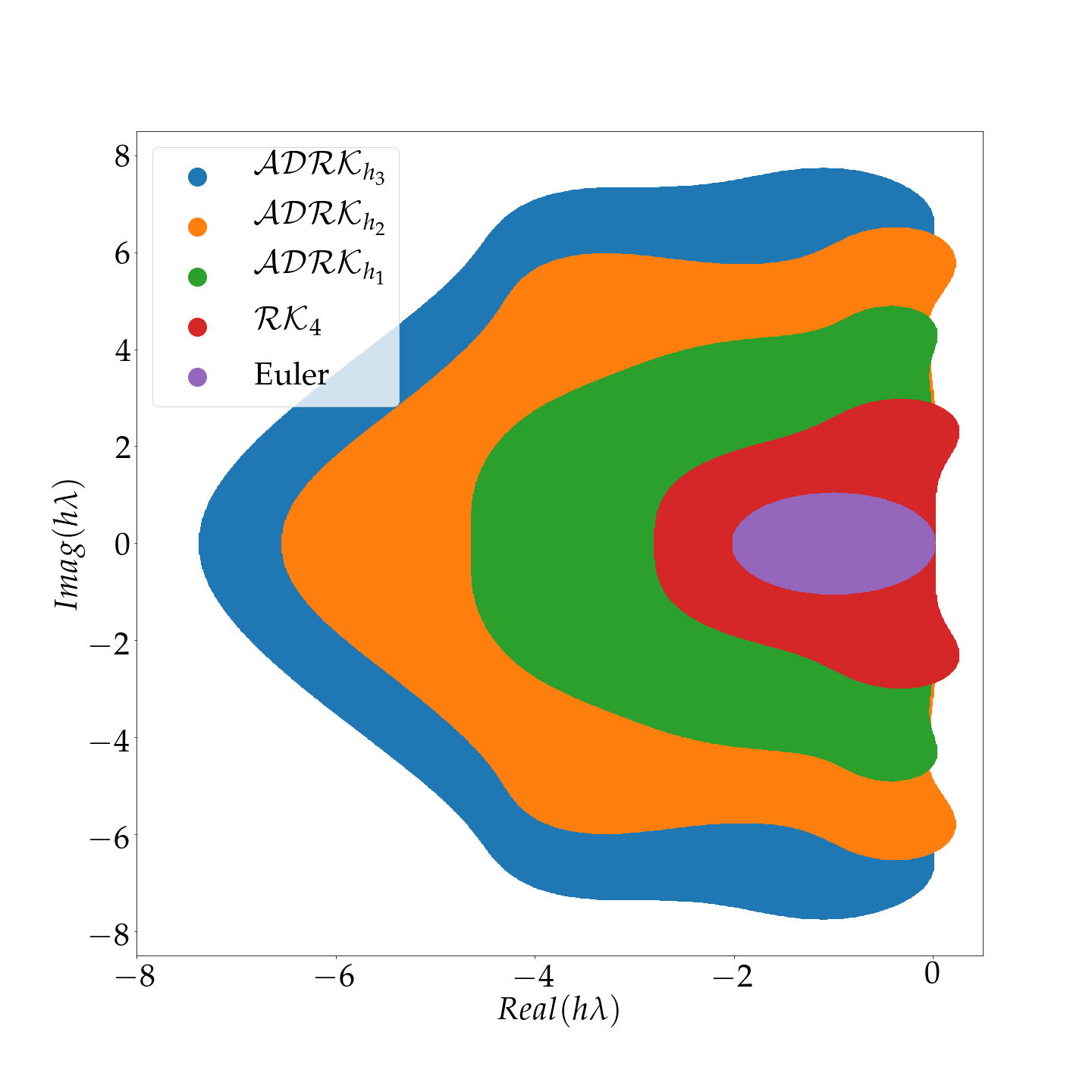}%
\caption{{{{ \bf  \em Stability region of the ADRK schemes (Lorenz 63 identification case-study)}.}
}}
\label{fig:C3_Stability_RINN_L63_dtall}
\end{figure}

\textbf{Stability and precision analysis of the learnt ADRK schemes}: Figure \ref{fig:C3_Stability_RINN_L63_dtall} illustrates the stability region of the trained schemes. The Euler and Runge-Kutta-4 schemes are shown as references. Interestingly, the higher the integration time-step of the trained model, the larger the stability region of the optimized scheme. Furthermore, analysing the stability polynomials of the ADRK schemes given in (\ref{EQ:C3_Taylor_gain_RINN_exp_L63}) reveals an interesting aspect of the proposed framework. The order of the integration scheme increases when moving from the data sampled at $h = 0.2$ to $h = 0.3$ but decreases when trained on the trajectory sampled at $h = 0.4$. These observations can be highlighted through the computation of the normalized coefficients error\footnote{This score is computed as the sum of the normalized root squared error of each coefficient of the stability polynomials, given in \ref{EQ:C3_Taylor_gain_RINN_exp_L63}, with respect to the true Taylor expansion of the analytical solution up to a given order $p$.} given in figure \ref{fig:c3_coeff_errors}. We believe that the order of the $\mathcal{ADRK}_{h_3}$ drops to lower the truncation error of the integration scheme as when considering high values of $\lambda h$, and as illustrated in figure \ref{fig:C3_Truncation_Error_Lin}, low order schemes tend to have a smaller truncation error, which guarantees a larger stability region, necessary for both the identification and the integration of the data driven model at $h = 0.4$. From this point of view, the learnt integration scheme was able to adapt, similarly to the previous case study to the dynamics of the observations, but also to their sampling leading to a correct identification of the dynamics. 


\begin{equation}
\label{EQ:C3_Taylor_gain_RINN_exp_L63}
\begin{split}
\centering
\mathcal{R}_{\mathcal{ADRK}_{h1}}(h\lambda)
              &= 1 + h\lambda + 0.5009 (h\lambda)^{2} + 0.1646 (h\lambda)^{3} + 0.04010 (h\lambda)^{4} + 0.007592 (h\lambda)^{5} +\\
              &0.001136 (h\lambda)^{6} + ...\\
\mathcal{R}_{\mathcal{ADRK}_{h2}}(h\lambda)
              &= 1 + h\lambda + 0.5021 (h\lambda)^{2} + 0.1691 (h\lambda)^{3} + 0.04273 (h\lambda)^{4} + 0.008826 (h\lambda)^{5} +\\
              &0.001510 (h\lambda)^{6} + 2.147^{-4}(h\lambda)^{7} + 2.550E^{-5}(h\lambda)^{8}+...\\
\mathcal{R}_{\mathcal{ADRK}_{h3}}(h\lambda)
              &= 1 + h\lambda + 0.4971 (h\lambda)^{2} + 0.1630 (h\lambda)^{3} + 0.03935 (h\lambda)^{4} + 0.007438 (h\lambda)^{5} +\\
              &0.001138 (h\lambda)^{6} + 1.411^{-4}(h\lambda)^{7} + 1.423E^{-5}(h\lambda)^{8}+...\\
\mathcal{R}_{exp}(h\lambda)
              &= 1 + h\lambda + 0.5 (h\lambda)^{2} + 0.16666 (h\lambda)^{3} + 0.04166 (h\lambda)^{4}+ 0.008333 (h\lambda)^{5} +\\
              &0.001388 (h\lambda)^{6} + 1.984^{-4}(h\lambda)^{7} + 2.480E^{-5}(h\lambda)^{8}+...
\end{split}
\end{equation}

\section{Conclusion}
\label{sec:3Conclusion}

In this work, we address the data-driven identification of numerical integration schemes for the simulation and the identification of ODE representations of dynamical systems. Formally, we focus on explicit Runge-Kutta integration schemes and we state the identification of these schemes as an optimization problem according to some predefined criteria. We implement numerically the proposed approach using deep learning framework to benefit from embedded automatic differentiation tools. Interestingly, we can consider both theory-guided and data-driven criterion. We demonstrate the relevance of the proposed framework on three case studies, namely the derivation of stable integration schemes for the time integration of linear equations, the identification of an integration scheme optimized for the resolution of a predefine non-linear ODE and the joint identification of numerical integration schemes and ODE representations from a sequence of observations.

The reported experiments demonstrate the relevance of the proposed framework. Specifically, we demonstrate experimentally that stability-constrained ADRK schemes can reach stability upper bounds found in related works \cite{van1972explicit,kinnmark1984one,ketcheson2013optimal} without resorting to any extra-parameterization or requiring to chose basis functions that generate stability functions as in \cite{ketcheson2013optimal}. We also show that trainable ADRK schemes can adapt to given dynamics and to the sampling of the observations to properly integrate non-linear dynamics. Besides, they make tractable the joint identification of an integration scheme and of an ODE representation of some observed dynamics. We may point out that the reported results do not relate to the parametric form of the ODE but truly to the fact that using classical integration schemes might not be appropriate for the considered identification task.


Future work will first investigate the proposed framework for the identification of data-driven representations of geophysical dynamics. Applications to satellite-derived sea surface observations appear very appealing. The typical sampling of these data range from a few hours to 10 days. The proposed ADRK then seem more appropriate compared with fixed step-size techniques to explore data-driven representations. 


From a methodological point of view, future work may first investigate the identification of ADRK schemes using other data-driven or stability-based criterions, including including ones defined as weighted versions of different criterions. Whereas this work focuses on explicit integration schemes, the proposed framework could be extended to implicit ones  including a non-linear differentiable solver. 
This could apply for instance to classical implicit solvers such as modified newton iteration \cite{cooper1993some}. 

\bibliographystyle{unsrtnat}
\bibliography{references}  






\section*{Acknowledgement}
This work was supported by Labex Cominlabs (grant SEACS), CNES (grant OSTST-MANATEE), Microsoft (AI EU Ocean awards), ANR Projects Melody and OceaniX. It benefited from HPC and GPU resources from Azure (Microsoft EU Ocean awards) and GENCI-IDRIS (Grant 2020-101030).

\clearpage
\appendix
\section{Integration schemes parameters}
\label{sec:integration_matrices}
We present in this section the integration coefficients {\em i.e.} the matrix $a$ and the vectors $b$ and $c$ of the trained integration schemes of section \ref{sec:identification_L63}.  Classically, Runge-Kutta coefficients are arranged in a Butcher table. The general form of this table is given in table \ref{tab:butcherTab}.
\begin{table}[h!]
\centering
\begin{tabular}{c|c}
     c   & a\\
     \hline
      & b
\end{tabular}
\caption {{\footnotesize{ \bf  \em Butcher table}.}}
\label{tab:butcherTab}
\end{table}

As an illustrative example, table \ref{tab:butcherTabRK45} resumes the well known Runge-Kutta 4 method.
\begin{table}[h!]
\centering
\begin{tabular}{c|cccc}
     0   & 0   & 0 & 0&0\\
     1/2 & 1/2 & 0 & 0&0\\
     1/2 & 0   & 1/2 & 0&0\\
     1   &0    & 0 & 1 &0\\
     \hline
      & 1/6 & 1/3 & 1/3 & 1/6
\end{tabular}
\caption {{\footnotesize{ \bf  \em Runge-Kutta 4 parameters}.}}
\label{tab:butcherTabRK45}
\end{table}

Tables \ref{tab:butcherTabRINNh1}, \ref{tab:butcherTabRINNh2} and \ref{tab:butcherTabRINNh3} highlight the parameters of the integration schemes $\mathcal{ADRK}_{h_1}$, $\mathcal{ADRK}_{h_2}$ and $\mathcal{ADRK}_{h_3}$ respectively trained on the Lorenz 63 data sampled from $h=0.2$ to $h=0.4$. 
\begin{table*}[h!]
\centering
\begin{adjustbox}{max width=\textwidth}
\begin{tabular}{c|ccccccccccc}
  0&0 &  0 &  0 &  0 &  0 &  0 &  0 & 0 & 0 & 0 & 0 \\
  0.00295&0.00295 &  0 &  0 &  0 &  0 &  0 &  0 & 0 & 0 & 0 & 0 \\
  0.00299&2.17592 & -2.17293 &  0 &  0 &  0 &  0 &  0 & 0 & 0 & 0 & 0 \\
  0.00300&2.98541 & -2.17397 & -0.80844 &  0 &  0 &  0 &  0 & 0 & 0 & 0 & 0 \\
  0.24293&-1.48213 &  3.52131 & -0.92958 & -0.86667 &  0 &  0 &  0 & 0 & 0 & 0 & 0 \\
  0.24293&-1.20813 &  0.81761 &  0.44993 &  0.09160 &  0.09192 &  0 &  0 & 0 & 0 & 0 & 0 \\
  0.32598&0.68880 & -0.37282 & -0.24275 & -0.02117 &  0.06104 &  0.21287 &  0 & 0 & 0 & 0 & 0 \\
  0.49628&0.28836 &  0.04642 & -0.17632 & -0.01809 & -0.01089 &  0.07776 &  0.28904 & 0 & 0 & 0 & 0 \\
  0.71550&0.05976 & -0.02908 &  0.45015 & -0.16537 & -0.08716 & -0.16371 &  0.10869 & 0.54223 & 0 & 0 & 0 \\
  0.83650&0.11298 &  0.23727 & -0.32220 & -0.10863 &  0.21937 &  0.49812 & -0.42969 & 0.34843 & 0.28086 & 0 & 0 \\
  1 &-0.52891 &  0.28477 &  0.30886 & -0.10163 &  0.19143 &  0.33904 &  0.06679 & 0.08134 & 0.14083 & 0.21750 & 0 \\
\hline
    &0.00013 & 0.00013 & 0.09499 & 0.00011 & 0.00033 & 0.03125 & 0.39447 & 0.06335 & 0.27661 & 0.02192 & 0.11670
\end{tabular}
\end{adjustbox}
\caption {{\footnotesize{ \bf  \em $\mathcal{ADRK}_{h_1}$ parameters}.}}
\label{tab:butcherTabRINNh1}
\end{table*}

\begin{table*}
\centering
\begin{adjustbox}{max width=\textwidth}
\begin{tabular}{c|ccccccccccc}
 0 &0 &  0 &  0 &  0 &  0 & 0 &  0 &  0 & 0 & 0 & 0 \\
 0.09168 &0.09168 &  0 &  0 &  0 &  0 & 0 &  0 &  0 & 0 & 0 & 0 \\
 0.09227 &-0.00598 &  0.09825 &  0 &  0 &  0 & 0 &  0 &  0 & 0 & 0 & 0 \\
 0.16187 &-0.23044 &  0.00477 &  0.38754 &  0 &  0 & 0 &  0 &  0 & 0 & 0 & 0 \\
 0.16204 &-0.25781 & -0.20551 &  0.79426 & -0.16890 &  0 & 0 &  0 &  0 & 0 & 0 & 0 \\
 0.43033  &-0.24147 &  0.53237 & -0.20302 &  0.46703 & -0.12458 & 0 &  0 &  0 & 0 & 0 & 0 \\
 0.43142 &1.16600 & -0.60220 & -0.35233 & -0.56929 &  0.27263 & 0.51660 &  0 &  0 & 0 & 0 & 0 \\
 0.60603 &0.34099 & -0.52091 &  0.27404 & -0.25958 &  0.25545 & 0.38264 &  0.13340 &  0 & 0 & 0 & 0 \\
 0.64954 &0.16978 & -0.07287 &  0.06919 & -0.17928 &  0.20630 & 0.26679 &  0.03296 &  0.15667 & 0 & 0 & 0 \\
 0.84016 &0.36515 & -0.36786 &  0.47056 & -0.52580 &  0.03298 & 0.45605 & -0.08132 & -0.12110 & 0.61150 & 0 & 0 \\
 1  &-0.11475 &  0.36853 &  0.01243 &  0.06163 & -0.10538 & 0.43067 & -0.02148 &  0.07745 & 0.14748 & 0.14341 & 0 \\
\hline
  &0.00010 & 0.00007 & 0.22873 & 0.00013 & 0.00020 & 0.27106 & 0.00003 & 0.00011 & 0.36341 & 0.04983 & 0.08633 \\
\end{tabular}
\end{adjustbox}
\caption {{\footnotesize{ \bf  \em $\mathcal{ADRK}_{h_2}$ parameters}.}}
\label{tab:butcherTabRINNh2}
\end{table*}

\begin{table*}
\centering
\begin{adjustbox}{max width=\textwidth}
\begin{tabular}{c|ccccccccccc}
  0&0 & 0 & 0 & 0 & 0 & 0 & 0 & 0 &  0 & 0 & 0 \\
  0.13282   &0.13282 & 0 & 0 & 0 & 0 & 0 & 0 & 0 &  0 & 0 & 0 \\
  0.18403&-0.01533 & 0.19936 & 0 & 0 & 0 & 0 & 0 & 0 &  0 & 0 & 0 \\
  0.28433&0.06830 & 0.03882 & 0.17721 & 0 & 0 & 0 & 0 & 0 &  0 & 0 & 0 \\
  0.37687&0.05088 & 0.04762 & 0.07653 & 0.20184 & 0 & 0 & 0 & 0 &  0 & 0 & 0 \\
  0.49570&0.05444 & 0.09223 & 0.08569 & 0.06398 & 0.19936 & 0 & 0 & 0 &  0 & 0 & 0 \\
  0.65230&0.10848 & 0.12104 & 0.08825 & 0.09688 & 0.01685 & 0.22079 & 0 & 0 &  0 & 0 & 0 \\
  0.66402&-0.00214 & 0.11741 & 0.10766 & 0.07588 & 0.12490 & 0.06258 & 0.17773 & 0 &  0 & 0 & 0 \\
  0.78754&0.03866 & 0.07149 & 0.12547 & 0.10885 & 0.07602 & 0.10123 & 0.06464 & 0.20117 &  0 & 0 & 0 \\
  0.85205&0.05165 & 0.10963 & 0.06966 & 0.11010 & 0.09470 & 0.14948 & 0.07204 & 0.10591 &  0.08888 & 0 & 0 \\
  1&0.02863 & 0.11511 & 0.05827 & 0.16788 & 0.08683 & 0.15047 & 0.02359 & 0.10025 & -0.09581 & 0.36479 & 0 \\
\hline
 &0.04237 & 0.10164 & 0.08155 & 0.12025 & 0.09046 & 0.13617 & 0.06988 & 0.11978 & 0.03882 & 0.14815 & 0.05093 \\
\end{tabular}
\end{adjustbox}
\caption {{\footnotesize{ \bf  \em $\mathcal{ADRK}_{h_3}$ parameters}.}}
\label{tab:butcherTabRINNh3}
\end{table*}

\section{Stability-constrained ADRK schemes, additional experiments}
\label{app:Stability_poly}
\begin{figure*}
\centering
  \subfloat[7-stage ADRK, Real axis inclusion]{%
\includegraphics[width=0.5\textwidth]{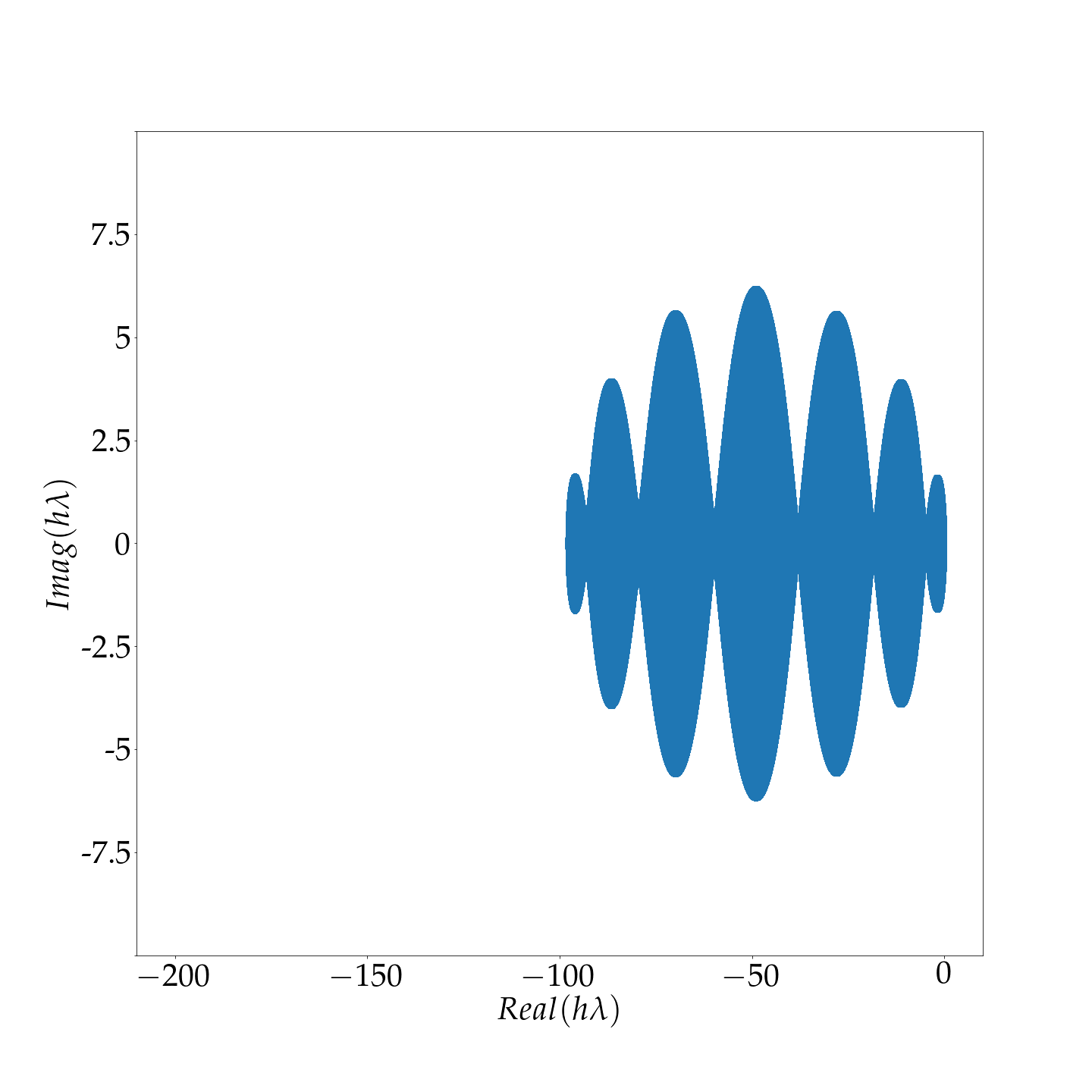}}
  \subfloat[7-stage ADRK, Imaginary axis inclusion]{%
\includegraphics[width=0.5\textwidth]{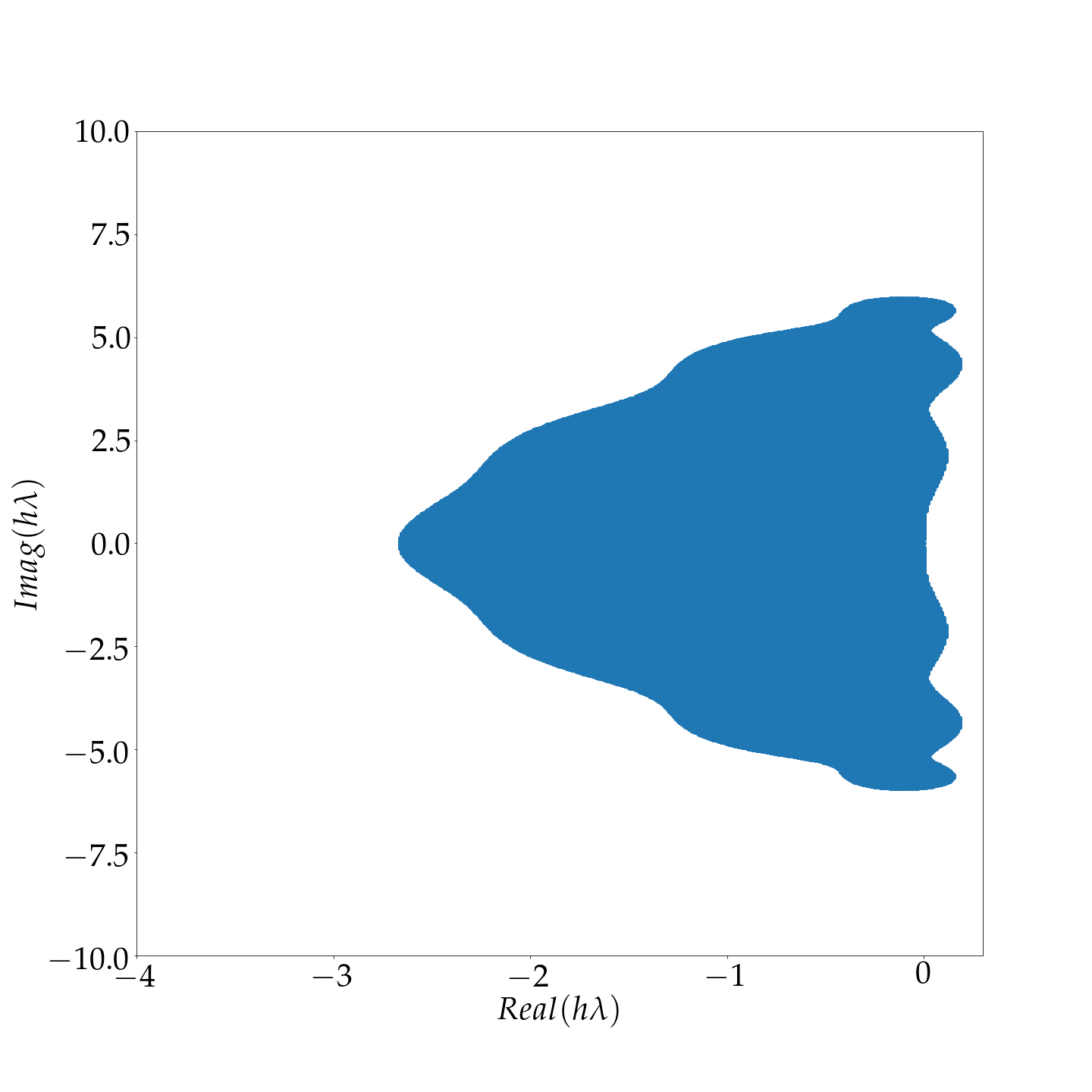}}\\
  \subfloat[10-stage ADRK, Real axis inclusion]{%
\includegraphics[width=0.5\textwidth]{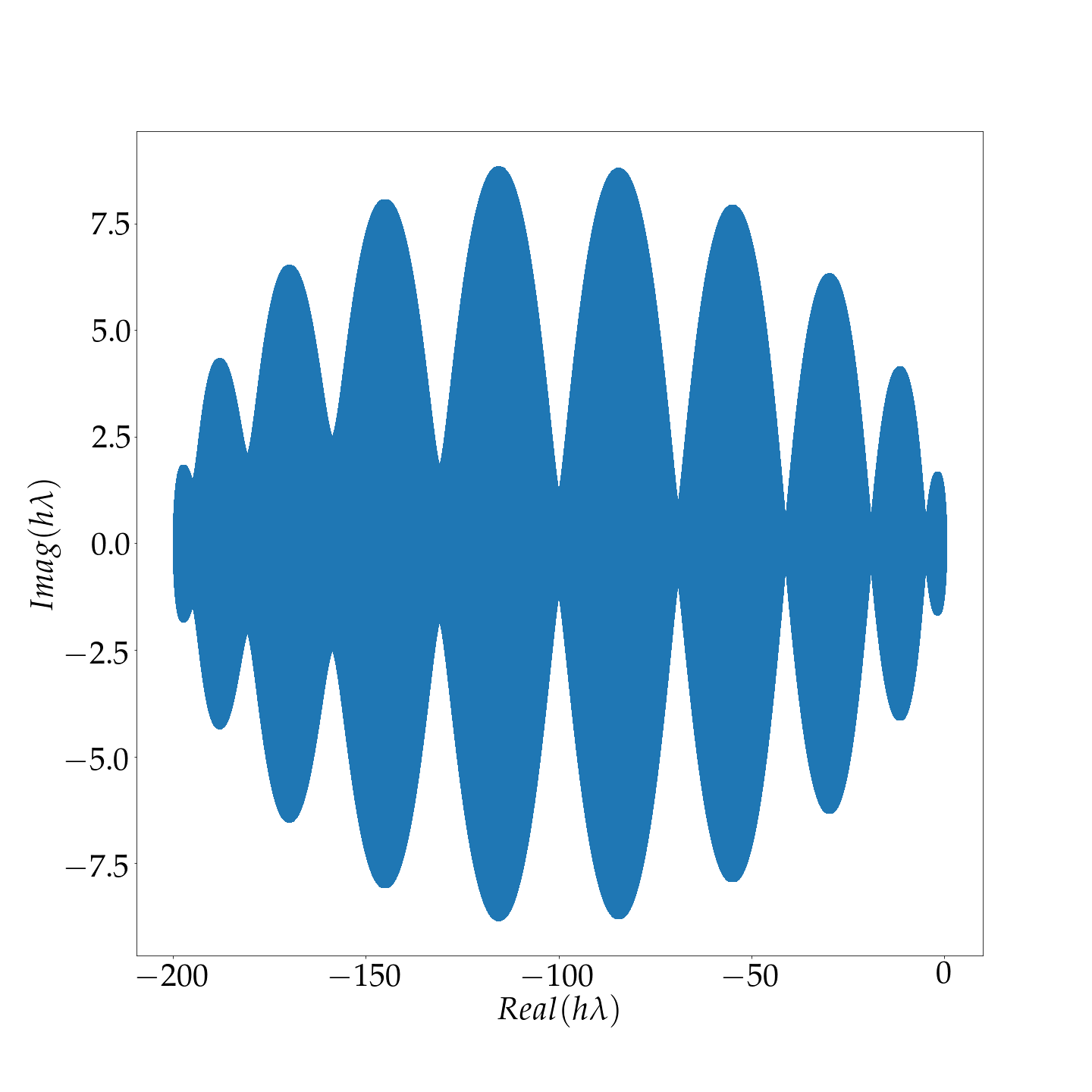}}
  \subfloat[10-stage ADRK, Imaginary axis inclusion]{%
\includegraphics[width=0.5\textwidth]{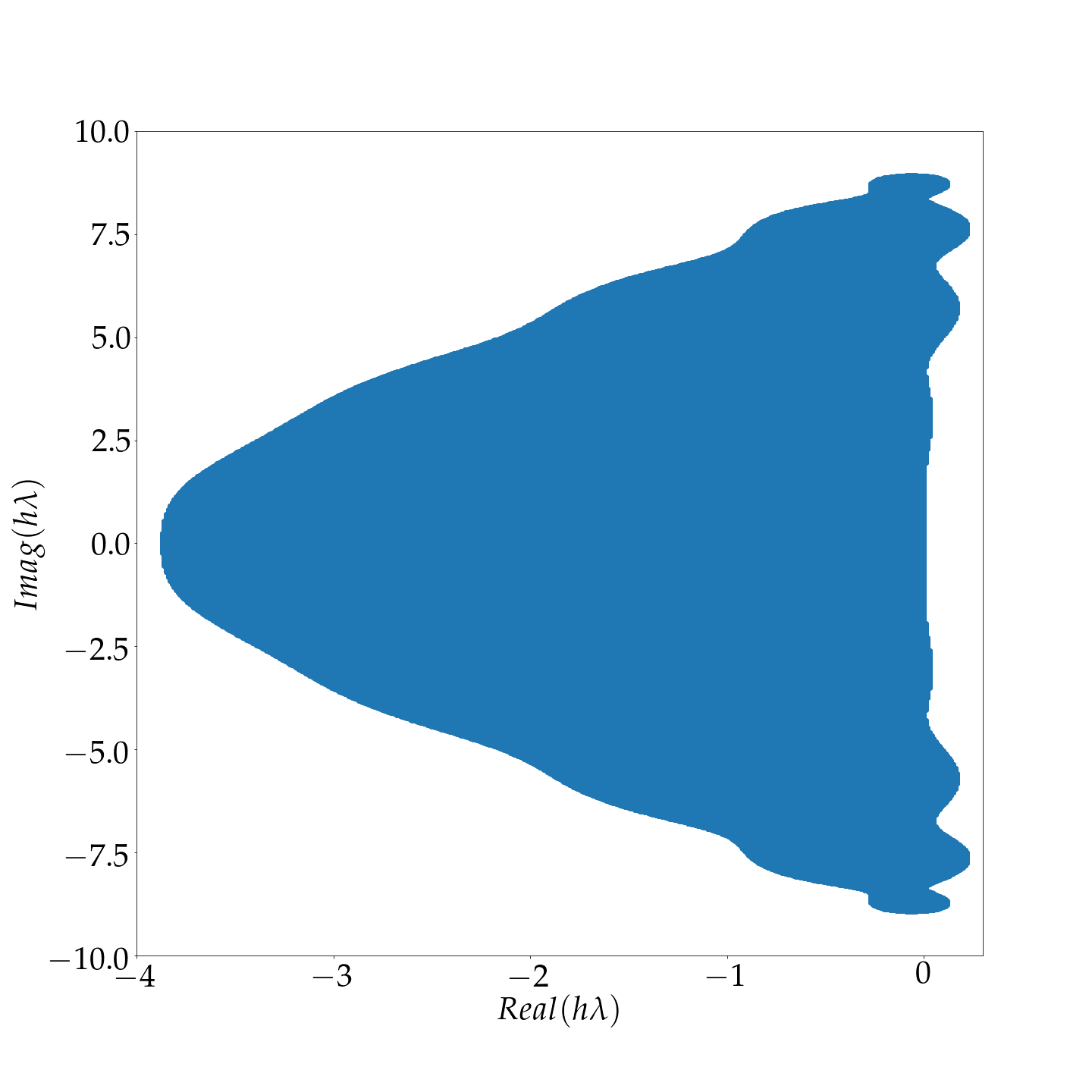}}
\caption{{{{ \bf  \em Two dimensional stability region of the stability-constrained $\mathcal{ADRK}$ integration schemes for several number of stages}: The tested, stability-constrained, $\mathcal{ADRK}$ schemes reach the stability limit of both real and imaginary axis inclusion case-studies.}}}
\label{fig:C3_stability_region_2d_Imag_ADRKb}
\end{figure*}

We complement the stability-constrained ADRK schemes experiment with two similar case studies based on an $\mathcal{ADRK}_{7}$ and an $\mathcal{ADRK}_{10}$. Specifically, both the Runge-Kutta schemes are optimized to reach the maximum possible value of $\lambda h = 2q^2$ for the real axis inclusion case study Figure \ref{fig:C3_stability_region_2d_Imag_ADRKb}.

\end{document}